\documentclass{amsart}

\usepackage[hidelinks]{hyperref}
\usepackage{color}
\usepackage{amssymb,amsmath,mathtools}
\usepackage{scalerel}
\usepackage{microtype}
\usepackage{lmodern}

\numberwithin{equation}{section}

\newtheorem{Thm}{Theorem}[section]
\newtheorem{Prop}[Thm]{Proposition}
\newtheorem{Lem}[Thm]{Lemma}
\newtheorem{Cor}[Thm]{Corollary}
\theoremstyle{definition}

\newtheorem{Rem}[Thm]{Remark}

\newcommand{\vol}{\operatorname{vol}}

\newcommand{\spa}{\operatorname{span}}

\newcommand{\id}{\mathrm{id}}

\renewcommand{\Re}{\operatorname{Re}}
\renewcommand{\Im}{\operatorname{Im}}
\newcommand{\Cl}{\operatorname{Cl}}
\newcommand{\End}{\operatorname{End}}
\newcommand{\Alt}{\operatorname{Alt}}
\newcommand{\sgn}{\operatorname{sgn}}

\newcommand{\wedges}{\wedge \cdots \wedge}
\newcommand{\bone}{\mathbf{1}}
\newcommand{\be}{\mathbf{e}}
\newcommand{\bi}{\mathbf{i}}
\newcommand{\bj}{\mathbf{j}}
\newcommand{\bk}{\mathbf{k}}
\newcommand{\bA}{\mathbf{A}}
\newcommand{\bG}{\mathbf{G}}
\newcommand{\bH}{\mathbf{H}}

\newcommand{\bO}{\mathbf{O}}
\newcommand{\bP}{\mathbf{P}}
\newcommand{\bS}{\mathbf{S}}

\newcommand{\C}{\mathbf{C}}
\newcommand{\N}{\mathbf{N}}
\newcommand{\R}{\mathbf{R}}

\newcommand{\ext}{\scalebox{1.0}[1.12]{\boldmath$\wedge$}}

\newcommand{\scirc}{\mathbin{\mkern-2mu\circ\mkern-2mu}}

\makeatletter
\newcommand{\extp}{\@ifnextchar^\@extp{\@extp^{\,}}}
\def\@extp^#1{\mathop{\bigwedge\nolimits^{\!#1}}}
\makeatother

\title{A calibration on $\mathbf R^{16}$ and Federer's product question}

\author{Roger Züst}

\email{roger\_zuest@hotmail.com}

\begin{document}

\begin{abstract}
Building upon the construction of a Cayley calibration adapted to a complex structure, we introduce a calibration $\Phi$ in $\ext^8 \mathbf R^{16}$ with $|\Phi^2| = 294$. This enables us to show that the product of two orthogonally supported calibrations is not necessarily a calibration, thereby providing a negative answer to a question posed by Federer. Dadok and Harvey developed a general method for constructing calibrations as outer products of two unit spinors in the Clifford algebra. We show that $\Phi$ arises from the product of two spinors with norm $1$ and $\sqrt{2}$.
\end{abstract}

 
\maketitle

%%%%%%%%%%%%%%%%%%%%%%%%%%%%%%%%%%%%%%%%%%%%%%%%%%%%%%%%%%%%%%%%%%%%%%%%%%%%%%%

\section{Introduction}

As introduced by Harvey and Lawson in \cite{HL}, a calibration on an $n$-dimensional Riemannian manifold $M$ is a closed $k$-form $\phi \in \Omega^k(M)$ that satisfies the pointwise inequality
\begin{equation}
	\label{eq:manifold_calibration}
\varphi_x(v_1 \wedges v_k) \leq \vol(v_1 \wedges v_k)
\end{equation}
for all $x \in M$ and all $v_1, \dots, v_k \in T_xM$. Calibrations are primarily used to identify volume-minimizing oriented submanifolds or, more generally, currents, whether in given homology classes or as solutions to the Plateau problem with prescribed boundary. An oriented $k$-dimensional submanifold $S \subset M$ is calibrated by $\varphi$, if $\varphi_x(v_1 \wedges v_k) = 1$ for all $x \in S$ and an oriented orthonormal basis $v_1, \dots, v_k \in T_x S$. In the terminology of Federer and Fleming \cite{FF}, inequality \eqref{eq:manifold_calibration} expresses the pointwise comass bound $\|\varphi_x\| \leq 1$, where $\varphi_x \in \ext^k T_xM$.

Let $X$ be a complex $2n$-dimensional Hermitian product space with holomorphic volume form $\Omega$ and Kähler form $\omega$. Define $\phi_{2n}$ as the sum of the special Lagrangian and Kähler calibrations:
\[
\phi_{2n} \coloneqq \Re\Omega + \tfrac{1}{n!}\omega^n \in \ext^{2n} X \cong \ext^{2n} \R^{4n}.
\]
For $n \geq 2$, $\phi_{2n}$ is itself a calibration. Specifically, for $n = 2$, this corresponds to a Cayley calibration in $\ext^{4} \R^{8}$, as established in \cite{HL}. For $n = 4$, we propose the following proper extension. The construction is based on the decomposition $X = X_1 \oplus X_2$ into complex vector spaces, each of dimension $\dim_\C(X_i) = 4$, and equipped with holomorphic volume forms $\Omega_1$ and $\Omega_2$, respectively.

\begin{Thm}
	\label{thm:main}
The form
\[
\Phi \coloneqq \Re\Omega + \tfrac{1}{4!}\omega^4 + \Re(\Omega_1 + \Omega_2) \wedge \tfrac{1}{2}\omega^2 \in \ext^8 X \cong \ext^8 \R^{16}
\]
is a self-dual calibration. Specifically, for all $\xi \in \bG(X,8)$, it satisfies
\[
|\Phi(\xi)| \leq |\Omega(\xi)| + \tfrac{1}{4!}|\omega^4(\xi)| + |(\Omega_1 \wedge \tfrac{1}{2}\omega^2)(\xi)| + |(\Omega_2 \wedge \tfrac{1}{2}\omega^2)(\xi)| \leq 1.
\]
Moreover, $\Phi^2 = 294 \vol$, where $\vol = \tfrac{1}{8!}\omega^8$ is the volume form on $X$.
\end{Thm}

The comass estimate above also implies that the specific choices of signs and phases for the components of $\Phi$ are not essential for it to define a calibration.

Federer raised the question of whether the equality $\|\alpha \wedge \beta\| = \|\alpha\| \|\beta\|$ holds when $\alpha$ and $\beta$ are supported in orthogonal subspaces \cite[\S1.8.4]{F}. Strictly speaking, the question in \cite{F} pertains to the mass, the dual norm of the comass. However, it remained open for the comass as well, as noted in \cite{M,DH} and \cite[Problem~10]{OPL}. Federer confirmed the validity of this equality when $\alpha$ or $\beta$ is simple. Additional positive results in low dimensions and low codimensions were given by Morgan in \cite{M}.

\begin{Thm}[{\cite[Theorem~5.1]{M}}]
	\label{thm:morgan}
Consider forms $\alpha \in \ext^k \R^m$, $\beta \in \ext^l\R^n$ for $0 \leq k \leq m$, $0 \leq l \leq n$. If $k \leq 2$, $m-k \leq 2$, $k=l=3$ or $m-k=n-k=3$, then the product $\alpha \wedge \beta \in \ext^{k+l}\R^{m+n}$ satisfies the equality $\|\alpha \wedge \beta\| = \|\alpha\| \|\beta\|$.
\end{Thm} 

The following proposition establishes a lower bound on the comass of certain wedge products. Its elementary proof is given in Section~\ref{subsec:exterior_algebras}.

\begin{Prop}
	\label{prop:cartesian}
Let $\alpha \in \ext^n \R^{2n}$ and define $\alpha_1,\alpha_2 \in \ext^{n}\R^{4n}$ as the isomorphic copies of $\alpha$ under the natural identifications $\R^{2n} \times \{0\} \cong \R^{2n} \cong \{0\} \times \R^{2n}$. Then, the comass of $\alpha_1 \wedge \alpha_2 \in \ext^{2n} \R^{4n}$ is bounded from below by
\[
\|\alpha_1 \wedge \alpha_2\| \geq \frac{|\alpha^2|}{2^{n}},
\]
where $\alpha^2 = \alpha \wedge \alpha \in \ext^{2n}\R^{2n}$ is a top form.
\end{Prop}

A sufficient condition for obtaining a counterexample to Federer's question is therefore
\begin{equation}
	\label{eq:main_estimate}
\frac{|\alpha^2|}{\|\alpha\|^2} > 2^n.
\end{equation}
This is satisfied by $\Phi$ from Theorem~\ref{thm:main}. Applying Proposition~\ref{prop:cartesian} yields the following conclusion.

\begin{Cor}
If $\Phi_1, \Phi_2 \in \ext^8 \R^{32}$ are orthogonally supported copies of $\Phi$, then their wedge product satisfies
\[
\|\Phi_1 \wedge \Phi_2\| \geq \frac{294}{256} > 1 = \|\Phi_1\| \|\Phi_2\|.
\]
\end{Cor}

%We want to emphasize, that this does not directly address the related questions concerning the mass, the dual norm of the comass, posed in \cite[\S1.8.4]{F}, nor the mass of Cartesian products of normal currents in \cite[\S4.1.8]{F}.

The term on the left-hand side of \eqref{eq:main_estimate} has applications in systolic geometry. In particular, the quantity $\mathrm{Wirt}_n$, defined as the maximum of $|\alpha^2|\|\alpha\|^{-2}$ over all nonzero elements $\alpha\in\ext^n \R^{2n}$, provides an upper bound for certain stable systolic ratios \cite{BKSW,HK,K}. Theorem~\ref{thm:main} implies that $\mathrm{Wirt}_8 \geq 294$. Further details are provided at the end of Section~\ref{sec:main}.

The $8$-planes calibrated by $\Phi$ are characterized in the following theorem. To simplify the formulation, we identify $X_1 = \C^4 = X_2$ and equip them with the standard holomorphic volume forms $\Omega_1$, $\Omega_2$ and Kähler forms $\omega_1$, $\omega_2$, respectively. The standard Kähler form on $\C^8$ is denoted by $\omega$. According to \cite[Lemma~6.13]{HL}, any $\xi \in \bG(\C^8,8)$ has a normal form representation as
\begin{align*}
\xi & = e_1 \wedge (\bi e_1 \cos \theta_1 + e_2 \sin \theta_1) \wedge e_3 \wedge (\bi e_3 \cos \theta_2 + e_4 \sin \theta_2) \; \wedge \\
 & \phantom{=} \; \; e_5 \wedge (\bi e_5 \cos \theta_3 + e_6 \sin \theta_3) \wedge e_7 \wedge (\bi e_7 \cos \theta_4 + e_8 \sin \theta_4) ,
\end{align*}
for a unitary basis $e_1, \dots, e_8$ of $\C^8$ and characteristic angles $0 \leq \theta_1 \leq \theta_2 \leq \theta_3 \leq \frac{\pi}{2}$, $\theta_3 \leq \theta_4 \leq \pi$. The complex $8 \times 8$ matrix whose columns are $e_1, \dots, e_8$ is denoted by $[e_1, \dots, e_8]$.

\begin{Thm}
	\label{thm:main2}
$\xi \in \bG(\C^8,8)$ is calibrated by
\[
\Phi = \Re(\Omega_1 \wedge \Omega_2) + \tfrac{1}{4!}\omega^4 + \Re(\Omega_1 + \Omega_2)\wedge \tfrac{1}{2}\omega^2
\]
if and only if one of the following four cases holds:
\begin{enumerate}
	\item $\xi$ is a special Lagrangian $8$-plane.
	\item $\xi$ is a complex $4$-plane.
	\item $\xi = \xi_1 \wedge \xi_2$, where each $\xi_i$ is a Cayley $4$-plane in $X_i$, i.e., a $4$-plane calibrated by $\Re \Omega_i + \tfrac{1}{2}\omega_i^2$.
	\item $\xi$ admits a normal form representation with $\theta_1 = \theta_2 = \theta_3 = \theta_4$ and
	\[
	[e_1, \dots, e_8] = \operatorname{diag}(U_1, U_2) S,
	\]
	where $S \in \mathrm{Sp}(4)$ and $U_1, U_2 \in \mathrm{SU}(4)$.
\end{enumerate}
\end{Thm}

We briefly explain why the compact symplectic group $\mathrm{Sp}(4)$ comes into play. If $r_1, \dots, r_8$ are the rows of $[e_1, \dots, e_8]$, then the matrix condition in (4) is equivalent to
\begin{equation*}
	\label{eq:main2}
\Re\tfrac{1}{2}\sigma^2 (r_1 \wedge r_2 \wedge r_3 \wedge r_4) = \Re\tfrac{1}{2}\sigma^2 (r_5 \wedge r_6 \wedge r_7 \wedge r_8) = 1 ,
\end{equation*}
where $\sigma = dz_1 \wedge dz_2 + \cdots + dz_7 \wedge dz_8$ is the standard $\C$-symplectic form on $\C^8$. Bryant and Harvey showed in \cite{BH} that $\Re\tfrac{1}{2}\sigma^2$ is a calibration and calibrated $4$-planes are contained in quaternionic $2$-planes.

Dadok and Harvey introduced a way to construct calibrations as products of unit spinors in \cite{DH}. Assume $(V,\langle\cdot,\cdot\rangle)$ is an oriented real inner product space of dimension $n = 8m$. As a vector space, the Clifford algebra $\Cl(V)$ is isomorphic to $\ext^\ast V$ and as an algebra, $\Cl(V)$ can be represented by $\End(\bP)$, where the pinor space $\bP$ is a real vector space of dimension  $16^m$. Up to a conformal factor, $\bP$ is equipped with a natural inner product $\varepsilon$ and has an orthogonal decomposition into positive and negative spinors $\bP = \bS^+ \!\oplus \bS^-$. These are the eigenspaces of the volume element $\nu \in \ext^n V$ with eigenvalues $\pm 1$. The outer product of two pinors $x,y \in \bP$ is $x \scirc y \in \End(\bP)$ defined by $(x \scirc y)(z) \coloneqq \varepsilon(y,z)x$ for all $z \in \bP$. The following statement is an extraction from results in \cite{DH}, primarily Proposition~2.1 and Theorem~2.2, which are pertinent to our discussion.

\begin{Thm}[\cite{DH}]
\label{thm:spinorprod_dh}
If $x,y \in \bS^+$ have unit norm, then
\[
\phi \coloneqq 16^m\, x \scirc y \in \End(\bP) \cong \Cl(V) \cong \ext_\ast V \cong \ext^\ast V
\]
satisfies:
\begin{enumerate}
	\item $\phi(\xi) = \varepsilon(x, \xi y) \leq |x||\xi y| \leq 1$ for all $k$ and all $\xi \in \bG(V,k) \subset \ext_k V$.
	\item $\phi = \sum_k \phi_{2k}$ for $\phi_{2k} \in \ext^{2k} V$ and each $\phi_{2k}$ is a calibration.
	\item $\phi = 1 + \cdots + \phi_{4k} + \cdots + \nu$ if $x = y$.
	\item $\ast \phi_{k} = \phi_{n-k}$ if $k \equiv 0 \text{ mod } 4$ and $\ast \phi_{k} = -\phi_{n-k}$ if $k \equiv 2 \text{ mod } 4$.
\end{enumerate}
\end{Thm}

Of particular interest in \cite{DH} are $\Cl(8)$ and $\Cl(16)$. In this work, $\Cl(8)$ is modelled by right-multiplications of elements in the octonions $\bO \cong \R^8$ acting on the pinor space $\bP(8) \coloneqq \bO^+ \!\oplus \bO^-$, where $\bO^\pm$ represent two copies of $\bO$. The usual convention is used that the standard basis $e_1, \dots, e_8 \in \R^8$ is such that $e_1$ is the unit element in $\bO$. The Clifford algebra $\Cl(16)$ can be represented as $\Cl(8) \otimes \Cl(8)$, acting on the pinor space $\bP(16) \coloneqq \bP(8) \otimes \bP(8)$. Our main observation in this context is the following statement.

\begin{Thm}
	\label{thm:spinorproduct}
Let $s \coloneqq e_1^+ \! \otimes e_1^+$, $s' \coloneqq e_2^+ \! \otimes e_2^+$ for the standard basis elements $e_1^+,e_2^+ \in \bO^+$ and define
\[
\phi \coloneqq 16^2\, (s + s')\scirc s \in \Cl(16) \cong \ext^\ast \R^{16} .
\]
Then $\phi = 1 + \cdots + \phi_{2k} + \cdots + \nu$ and every $\phi_{2k}$ is a calibration. Moreover, $\R^{16}$ admits a complex structure such that $\phi_8$ takes the form given in Theorem~\ref{thm:main}.
\end{Thm}

Squared spinors in the algebra $\Cl(16)$, such as $\psi = 16^2 s \scirc s$, along with their associated calibrated geometries are discussed in \cite[Theorem~4.1]{DH}. Specifically, $\psi$ corresponds to the $\mathrm{Spin}(7) \times \mathrm{Spin}(7)$ geometry. Using Theorem~\ref{thm:spinorproduct}, we obtain extensions of the calibrations $\psi_k$. Moreover, since $|s| = 1$ and $|s + s'| = \sqrt{2}$, this illustrates that some calibrations can be produced via spinor products without requiring unit norm.

%\medskip
%\noindent \textsc{Acknowledgment.} I would like to thank the referees for their thorough review and independent verification of Theorem 1.1.

\section{Preliminaries}

\subsection{Exterior algebras}
\label{subsec:exterior_algebras}

We primarily follow the notation of \cite[\S1]{F}, where further details can be found. Let $X$ be an $n$-dimensional real Euclidean vector space. The space of $m$-covectors $\ext^m X$ consists of all alternating $\R$-multilinear forms mapping $m$ vectors in $X$ to $\R$. For the forms with values in $\C$, the space is denoted by $\ext^m (X,\C)$. The space $\ext^m X$ is naturally identified with $(\ext_m X)^*$, the dual space of $m$-vectors. The pairing between $\alpha \in \ext^m X$ and $\xi \in \ext_m X$ is written as $\alpha(\xi) \in \R$. Both $\ext^m X$ and $\ext_m X$ have dimension $\binom{n}{m}$. An element $\xi \in \ext_m X$ is called simple if there are $x_1, \dots, x_m \in X$ with $\xi = x_1 \wedges x_m$. The inner product $\langle \cdot, \cdot \rangle$ on $X$ extends naturally to the full exterior algebras $\ext_\ast X = \bigoplus_{m=0}^n \ext_m X$ and $\ext^\ast X = \bigoplus_{m=0}^n \ext^m X$. These extensions are characterized by the orthogonality of components of different degree and the identity
\[
\langle x_1 \wedges x_m, y_1 \wedges y_m\rangle = \det(\langle x_i, y_j\rangle)_{1\leq i, j \leq m}.
\]
The associated Euclidean norms are denoted by $|\cdot|$. The Grassmannian of oriented $m$-planes in $X$ is
\[
\bG(X,m) \coloneqq \{\xi \in \ext_m X : \xi \text{ is simple and } |\xi| = 1 \} .
\]
The comass of $\alpha \in \ext^m X$ is defined by $\|\alpha\| \coloneqq \sup\bigl\{\alpha(\xi) : \xi \in \bG(X,m) \bigr\}$. Equivalently, when interpreting $\alpha$ as a multilinear form, its comass is
\begin{equation*}
\|\alpha\| = \sup\bigl\{\alpha(x_1, \dots, x_m) : |x_i| = 1 \text{ for all } i \bigr\} .
\end{equation*}

Let $T \colon X^m \to \R$ be a multilinear functional, i.e., a tensor in $\otimes^m X^\ast$. Its al\-ter\-na\-tion is defined by
\[
\Alt(T)(x_1, \dots, x_{m}) \coloneqq \frac{1}{m!}\sum_{\sigma \in S_{m}} \sgn(\sigma) T(x_{\sigma(1)}, \dots, x_{\sigma(m)}) .
\]
Given $\alpha \in \ext^{k}X$ and $\beta  \in \ext^{l}X$, their wedge product is defined by
\[
\alpha \wedge \beta \coloneqq \frac{(k+l)!}{k!l!} \Alt(\alpha \otimes \beta) \in \ext^{k+l}X.
\]
For an integer $m \in \N$, set $\N_m \coloneqq \{1, \dots, m\}$ and $\N_0 \coloneqq \emptyset$. If $0 \leq k \leq m$, then $\Lambda(m,k)$ denotes the set of all increasing sequences $\N_k \to \N_m$. If $k+l = m$ and $I \in \Lambda(m,k)$, then $I^c$ is the unique element of $\Lambda(m,l)$ with $I(\N_k) \cup I^c(\N_l) = \N_m$. The permutation $\sigma_I = (I,I^c) \in S_m$ is  obtained by concatenating $I$ and $I^c$. The wedge product, as defined above, can alternatively be expressed as
\[
(\alpha \wedge \beta)(x_1, \dots, x_m) = \sum_{I \in \Lambda(m,k)} \sgn(\sigma_I)\, \alpha(x_{I_1}, \dots, x_{I_{k}})\, \beta(x_{I_1^c}, \dots, x_{I_{l}^c}).
\]
This is merely a restatement of \cite[\S1.4.2]{F}, where shuffles are employed. The following more general result will also be required.

\begin{Lem}
	\label{lem:tensor_alternation}
Assume that $S \in \otimes^k X^\ast$, $T \in \otimes^l X^\ast$ and $m = k+l$. Define $U \in \otimes^m X^\ast$ by
\[
U(x_1, \dots, x_{m}) \coloneqq \sum_{I \in \Lambda(m,k)} \sgn(\sigma_I)\, S(x_{I_1}, \dots, x_{I_{k}})\, T(x_{I_1^c}, \dots, x_{I_{l}^c}) .
\]
Then
\[
\Alt(U) = \frac{m!}{k!l!} \Alt(S \otimes T) = \Alt(S) \wedge \Alt(T) .
\]
\end{Lem}

Since the argument reduces to a routine computation, we omit the details. The essential point is that each $\sigma \in S_m$ can be uniquely expressed as $\sigma_I \circ \sigma_A \circ \sigma_B$, where $\sigma_A$ and $\sigma_B$ are permutations of $A \coloneqq \N_k$ and $B \coloneqq \N_m \setminus A$, respectively, and $I \in \Lambda(m,k)$.

In case $X$ is additionally oriented with volume form $\vol \in \ext^n X$, the Hodge star operator $\ast \colon \ext^{k} X \to \ext^{n-k} X$ is defined by the relation
\[
\alpha \wedge \!\ast \mu = \langle \alpha, \mu \rangle \vol
\]
for all $\alpha, \mu \in \ext^k X$. If $\vol = E_1 \wedges E_n$, where $E_1, \dots, E_n \in X^\ast$ is an orthonormal basis, then
\begin{equation}
	\label{eq:hodge_basis}
	\ast E_I = \sgn(\sigma_I) E_{I^c}
\end{equation}
for all $I \in \Lambda(n,k)$. In case $\dim(X) = 2k$ is even, a form $\mu \in \ext^k X$ is called self-dual if $\ast \mu = \mu$.

Since the proof of Proposition~\ref{prop:cartesian} uses only the basic properties of the wedge product, we include it here.

\begin{proof}[Proof of Proposition~\ref{prop:cartesian}]
Let $\alpha_1, \alpha_2 \in \ext^{n}\R^{4n}$ denote the two orthogonal copies of $\alpha \in \ext^{n}\R^{2n}$, as given in the statement. Let $e_i \in \R^{2n}$, $f_i \in \R^{2n} \times \{0\}$ and $g_i \in \{0\} \times \R^{2n}$ be the $i$th standard basis elements. Set $v_i \coloneqq \tfrac{1}{\sqrt{2}}(f_i + g_i)$ and define
\[
\xi \coloneqq v_1 \wedges v_{2n} \in \bG(\R^{4n},2n) .
\]
From the definition of the wedge product, we obtain
\begin{align*}
(\alpha_1 \wedge \alpha_2)(\xi) & = \sum_{I \in \Lambda(2n,n)} \sgn(\sigma_I)\, \alpha_1(v_{I_1} \wedges v_{I_{n}})\, \alpha_2(v_{I_1^c} \wedges v_{I_{n}^c}) \\
& = \sum_{I \in \Lambda(2n,n)} \sgn(\sigma_I)\, \alpha_1 \bigl(\tfrac{1}{\sqrt{2}} f_{I_1} \wedges \tfrac{1}{\sqrt{2}}f_{I_{n}} \bigr)\, \alpha_2 \bigl( \tfrac{1}{\sqrt{2}} g_{I_1^c} \wedges \tfrac{1}{\sqrt{2}} g_{I_{n}^c} \bigr) \\
& = \frac{1}{2^{n}} \sum_{I \in \Lambda(2n,n)} \sgn(\sigma_I)\, \alpha(e_{I_1} \wedges e_{I_{n}})\, \alpha(e_{I_1^c} \wedges e_{I_{n}^c}) \\
& = \frac{1}{2^{n}} (\alpha \wedge \alpha)(e_1 \wedges e_{2n}) .
\end{align*}
Thus,
\[
\|\alpha_1 \wedge \alpha_2\| \geq |(\alpha_1 \wedge \alpha_2)(\xi)| \geq \frac{|\alpha^2|}{2^{n}} 
\]
as claimed.
\end{proof}

\subsection{Normed algebras}
\label{subsec:normed_algebras}

A normed algebra $\bA$ is a finite-dimensional Euclidean vector space over $\R$ endowed with a bilinear product and a multiplicative unit $\bone$, such that $|xy| = |x||y|$ holds for all $x,y \in \bA$. In particular, the unit satisfies $|\bone|=1$. Hurwitz's theorem asserts that every normed algebra is isomorphic to one of the following: the real numbers $\R$, the complex numbers $\C \cong \R^2$, the quaternions $\bH \cong \R^4$, or the octonions $\bO \cong \R^8$. The algebra $\bH$ is associative but not commutative, while $\bO$ is neither associative nor commutative. Each of these algebras admits an orthogonal decomposition $\bA = \Re \bA \oplus \Im \bA$, where $\Re \bA$ is spanned by $\bone$ and is identified with $\R$. Conjugation is defined by $\bar{x} \coloneqq \Re x - \Im x$. For further details, we refer to \cite[Appendix~IV.A]{HL}. It always holds
\begin{equation}
	\begin{gathered}
	\label{eq:algebra_basics}
	\bar{\bar{x}} = x, \quad \overline{xy} = \bar{y}\bar{x}, \quad \bar{x}x = |x|^2, \\
	\langle x,yz \rangle = \langle x\bar{z}, y\rangle, \quad \langle x,zy \rangle = \langle \bar{z}x, y\rangle.
	\end{gathered}
\end{equation}
If moreover $x$ and $y$ are orthogonal, then
\begin{equation}
	\label{eq:oct_permute}
x(\bar{y}z) = -y(\bar{x}z), \quad (z\bar{x})y = -(z\bar{y})x .
\end{equation}

We briefly introduce the notation and fundamental properties of these algebras, as well as vector spaces over them. For more details, see \cite{F, HL, H}.

As a vector space over $\R$, $\C^n$ is identified with $\R^{2n}$ through
\[
(x_1,x_2, \dots, x_{2n-1},x_{2n}) \mapsto (x_1 + \bi x_2, \dots, x_{2n-1} + \bi x_{2n}) .
\]
The Hermitian product of $v = (v_1, \dots, v_n)$ and $w = (w_1, \dots, w_n)$ in $\C^n$ is defined by
\[
H(v,w) \coloneqq \sum_{j=1}^n \bar{v}_j w_j \in \C.
\]
Let $Z_1, \dots, Z_n \in \ext^1(\C^n,\C)$ be the coordinate functions on $\C^n$, and $E_1, \dots, E_{2n} \in \ext^1 \R^{2n}$ the coordinate functions on $\R^{2n}$, with $Z_j = E_{2j-1} + \bi E_{2j}$. The Hermitian product can be decomposed into real and imaginary part as
\[
H(v,w) = \langle v, w \rangle + \bi \omega(v,w) ,
\]
where $\langle \cdot, \cdot \rangle$ is the usual inner product on $\R^{2n}$ and $\omega$ is the Kähler form
\[
\omega = \frac{\bi}{2} \sum_{j=1}^n Z_j \wedge \bar{Z}_j = \sum_{j=1}^n E_{2j-1} \wedge E_{2j} \in \ext^2 \C^n .
\]

The standard basis of the quaternions $\bH \cong \R^4$ is denoted by $\bone, \bi, \bj, \bk$, and satisfies the Hamiltonian identities $\bi^2 = \bj^2 = \bk^2 = \bi\bj\bk = -\bone$. The quaternionic Hermitian form on $\bH^n$ is given by
\[
\varepsilon(p,q) \coloneqq \sum_{l=1}^n \bar{p}_l q_l \in \bH .
\]
If $\C^{2n}$ is identified with $\bH^n$ via
\[
(z_1, \dots, z_{2n}) \mapsto (z_1 + \bj z_2, \dots,  z_{2n-1} + \bj z_{2n}) ,
\]
then $\varepsilon = H + \bj\sigma$, where
\[
\sigma = Z_1 \wedge Z_2 + \cdots + Z_{2n-1} \wedge Z_{2n} \in \ext^2 (\C^{2n},\C)
\]
is the standard $\C$-symplectic form on $\C^{2n}$. The compact symplectic group $\mathrm{Sp}(n)$ consists of those complex $2n \times 2n$ matrices that preserve $\varepsilon$, i.e., $A \in \mathrm{Sp}(n)$ if $\varepsilon(Ap,Aq) = \varepsilon(p,q)$. Thus, $\mathrm{Sp}(n)$ preserves both $H$ and $\sigma$. The group $\mathrm{Sp}(n)$ can be expressed as $\mathrm{Sp}(2n,\C) \cap \mathrm{U}(2n)$, where the $\C$-symplectic group $\mathrm{Sp}(2n,\C)$ consists of those complex $2n \times 2n$ matrices $S$ that satisfy $S^T J_{2n} S = J_{2n}$, with $J_{2n} \coloneqq \operatorname{diag}(J_2, \dots, J_2)$, where
\[
J_2 \coloneqq \begin{pmatrix}
	0 & 1 \\
	-1 & 0
\end{pmatrix} .
\]
By definition, $\mathrm{Sp}(2n,\C)$ is the group that fixes $\sigma$. for all $S \in \mathrm{Sp}(n)$, $\det(S) = 1$ holds, and thus $\mathrm{Sp}(n) = \mathrm{Sp}(2n,\C) \cap \mathrm{SU}(2n)$. We will also use that $\mathrm{Sp}(n)$ is closed under transposition, a fact that is readily available with the description of $\mathrm{Sp}(2n,\C)$ above.

If $\bO$ is identified with $\bH \oplus \bH \be$, where the standard orthonormal basis $e_1, \dots, e_8$ is given by $\bone$, $\bi$, $\bj$, $\bk$, $\be$, $\bi\be$, $\bj\be$, $\bk\be$, all the multiplication rules in $\bO$ follow from the Hamiltonian identities of $\bH$, $\bone x = x\bone = x$ and \eqref{eq:oct_permute}. In particular, for distinct $i,j \geq 2$, we have $e_i^2 = - e_1$, $e_ie_j = -e_je_i$ and if $x,y \in \{\bi, \bj, \bk\}$ are distinct, then
\begin{equation}
	\label{eq:octonion_rules}
(x\be)x = \be, \quad \be(x \be) = x, \quad x(y \be) = (yx)\be, \quad (x\be)(y \be) = yx.
\end{equation}

The following explicit computation in $\bO$ will be used.

\begin{Lem}
	\label{lem:octonion_product}
For $n \geq 0$, define $Q_n \colon \bO^n \to \bO$ recursively by setting $Q_0 \coloneqq \bone$ and
\[
Q_n(x_1, \dots, x_n) \coloneqq Q_{n-1}(x_2, \dots, x_n)x_1 .
\]
Then $Q_8(e_1, \dots, e_8) = \bone$.
\end{Lem}

\begin{proof}
With \eqref{eq:octonion_rules}, one readily obtains
\[
Q_4(\be, \bi\be, \bj\be, \bk\be) = Q_4(\bone, \bi, \bj, \bk) = \bone.
\]
Thus $Q_8(e_1, \dots, e_8) = \bone$.
\end{proof}

The threefold cross product on $\bO$ is defined by 
\[
x \times y \times z \coloneqq \frac{1}{2}(x(\bar{y}z) - z(\bar{x}y))
\]
in \cite[Definition~B.3]{HL} and shown to be alternating in \cite[Lemma~B.4]{HL}. The fourfold cross product on $\bO$ is defined by
\[
x \times y \times z \times w \coloneqq \tfrac{1}{4}(\bar{x}(y \times z \times w) - \bar{y}(z \times w \times x) + \bar{z}(w \times x \times y) - \bar{w}(x \times y \times z))
\]
in \cite[Definition~B.6]{HL} and is also shown to be alternating in \cite[Lemma~B.7]{HL}. If $x,y,z,w \in \bO$ are orthogonal, then by \cite[(B.5) and (B.8)]{HL},
\begin{equation}
	\label{eq:cross}
y \times z \times w = y(\bar{z}w) \quad \text{and} \quad x \times y \times z \times w = \bar{x}(y(\bar{z}w)) .
\end{equation}
The Cayley calibration $\phi_{\mathrm{Cay}} \in \ext^4 \R^8$ is defined by
\begin{equation}
	\label{eq:cayley_def}
\phi_{\mathrm{Cay}}(x,y,z,w) \coloneqq \langle x, y \times z \times w\rangle = \langle x \times y \times z \times w, \bone \rangle ,
\end{equation}
see \cite[Definition~1.21]{HL} and \cite[Lemma~B.9]{HL}.

\begin{Lem}
	\label{lem:cross_alternation}
$\phi_{\mathrm{Cay}} = \Alt(T)$, where $T(x_1,x_2,x_3,x_4) \coloneqq \langle \bar{x}_1(x_2(\bar{x}_3 x_4)), \bone \rangle$.
\end{Lem}

\begin{proof}
$\Alt(T)$ and $\phi_{\mathrm{Cay}}$ are equal when evaluated on orthogonal systems by \eqref{eq:cross} and \eqref{eq:cayley_def}. Since both forms are alternating, they are equal on any system.
\end{proof}

Expressed in the standard dual basis of $\R^8$,
\begin{align}
	\nonumber
\phi_{\mathrm{Cay}} & = +\, E_{1234} + E_{1256} - E_{1278} + E_{1357} + E_{1368} + E_{1458} - E_{1467} \\
	\label{eq:cayley_basis}
& \phantom{=}\ + E_{5678} + E_{3478} - E_{3456} + E_{2468} + E_{2457} + E_{2367} - E_{2358} ,
\end{align}
see \cite[Corollary~1.31]{HL}.

Note that the terminology \emph{a} Cayley calibration, in contrast to \emph{the} standard Cayley calibration ($\phi_{\mathrm{Cay}}$), is used to denote any self-dual calibration in $\ext^4 \R^8$ whose isotropy subgroup is isomorphic to $\operatorname{Spin}(7)$. It follows from \cite[Theorem~3.7]{DHM} that any such Cayley calibration lies in the $\operatorname{SO}(8)$-orbit of $\phi_{\mathrm{Cay}}$ or $-\phi_{\mathrm{Cay}}$.

\begin{Rem}
Starting with $\R$, the Cayley-Dickson process successively produces $\C$, $\bH$ and $\bO$ and further algebra structures in dimensions $2^n$. Based on $\bA$, the product on $\bA \oplus \bA$ is defined by
\[
(a,b)(c,d) \coloneqq (ac - \bar{d}b, da + b\bar{c}) .
\]
The Cayley calibration is derived through the alternation of an iterated octonionic product. A potential approach to obtaining calibrations in $\ext^8 \R^{16}$ could involve the alternation of products on the sedenions $\bS \cong \bO \oplus \bO \cong \R^{16}$, arising from the Cayley–Dickson construction. However, this algebra contains zero divisors and is not normed, which seems to obstruct the formation of forms with favorable properties constructed in this way.
\end{Rem}

\subsection{Classical calibrations}
\label{subsec:classical_calibrations}

We shortly review two primary examples of calibrations presented in \cite{HL}. Assume $Z_j \in \ext^1(\C^n,\C)$ and $E_j \in \ext^1 \R^{2n}$ are as above. The special Lagrangian calibration is the real part of the holomorphic volume form
\[
\Omega = Z_1 \wedges Z_n = (E_1 + \bi E_2) \wedges (E_{2n-1} + \bi E_{2n}) \in \ext^n(\C^n,\C) .
\]
The form $\Omega$ is actually multilinear over $\C$, so an element of $\ext^n_\C \C^n$ in the notation of \cite[\S1.6.6]{F}. The product Kähler forms are $\omega_k \coloneqq \tfrac{1}{k!}\omega^k \in \ext^{2k} \C^n$ for $k=0, \dots, n$. In the following lemma, $\R^{2n}$ has the standard orientation with volume form $\omega_n = E_1 \wedges E_{2n}$.

\begin{Lem}
	\label{lem:standard_calibrations}
The forms $\omega_k$ and $\Re \Omega$ are expressed as linear combinations of the basis elements $E_I$, $I \in \Lambda(2n,\cdot)$, with coefficients $\pm1$. Moreover,
\begin{enumerate}
	\item $|\omega_k|^2 = \binom{n}{k}$, and $\ast \omega_k = \omega_{n-k}$ for all $k = 0, \dots,  n$. In particular, $\omega_{n/2}$ is self-dual when $n$ is even.
	\item $|\!\Re \Omega|^2 = 2^{n-1}$ and $\Re \Omega$ is self-dual when $n$ is even.
\end{enumerate}
\end{Lem}

\begin{proof}
If $\Lambda_k \subset \Lambda(2n,2k)$ is the collection of all increasing sequences composed of pairs $(2j-1,2j)$, then
\[
\omega_k = \sum_{I \in \Lambda_k} E_{I} .
\]
Since $\Lambda_{k}$ corresponds to the set of all possible selections of $k$ pairs from the $n$ available pairs,
\[
|\omega_k|^2 = \binom{n}{k} .
\]
Because $\sgn(\sigma_I) = 1$ for all $I \in \Lambda_k$, it follows immediately from \eqref{eq:hodge_basis} that $\ast \omega_k = \omega_{n-k}$. 

For any $I \in \Lambda(2n,n)$, let $N(I)$ be the number of even indices that appear in the range of $I$. These indices correspond to the imaginary part of complex coordinates. Define $\Lambda(\Omega) $ as the set of all $I \in \Lambda(2n,n)$ satisfying the conditions that $N(I)$ is even and $|I \cap (2j-1,2j)| = 1$ for all $j=1, \dots, n$. Then
\[
\Re \Omega = \sum_{I \in \Lambda(\Omega)} (-1)^{N(I)/2} E_I .
\]
The signs that appear correspond to powers of $\bi$. Thus
\[
|\! \Re \Omega|^2 = |\Lambda(\Omega)| = \sum_{m \text{ \scriptsize even}} \binom{n}{m} = 2^{n-1} .
\]
The binomial coefficient $\binom{n}{m}$ represents the number of ways to choose $m$ even indices from the $n$ available pairs.

Clearly, $N(I) + N(I^c) = n$ for all $I \in \Lambda(2n,n)$. Thus $\Lambda(\Omega)$ is closed under complement when $n$ is even. Moreover, if $I \in \Lambda(\Omega)$ and $n$ is even, then
\[
(-1)^{N(I)/2}(-1)^{N(I^c)/2} = (-1)^{n/2} .
\]
For $I \in \Lambda(\Omega)$, we have $\sgn(\sigma_I) = (-1)^{l(I)}$, where
\[
l(I) = (n-1) + \delta_1 + (n-2) + \delta_2 + \cdots ,
\]
with $\delta_j = 0$ if $2j \in I^c$ and $\delta_j = 1$ if $2j-1 \in I^c$. Thus
\[
l(I) = \tfrac{1}{2}n(n-1) + \sum_i \delta_i .
\]
The sum $\sum_i \delta_i$ counts the number of odd indices in $I^c$, which is equal to the number of even indices in $I$. Thus $\sum_i \delta_i = N(I)$, which is even. Consequently, if $n$ is even, we have $\sgn(\sigma_I) = (-1)^{l(I)} = (-1)^{n/2}$. Using \eqref{eq:hodge_basis}, it follows that
\begin{align*}
\ast\Re \Omega & = \sum_{I \in \Lambda(\Omega)} (-1)^{N(I)/2} \ast\! E_I = \sum_{I \in \Lambda(\Omega)} (-1)^{N(I)/2} (-1)^{n/2} E_{I^c} \\
 & = \sum_{I \in \Lambda(\Omega)} (-1)^{N(I^c)/2} E_{I^c} = \Re \Omega .
\end{align*}
This concludes the proof.
\end{proof}

For the calculations in the final section, we require additional preparation concerning the Cayley form and its connection to certain complex structures. Let $J$ be the complex structure on $\bO$, defined by $Jv \coloneqq v \bi$ (note that $\bi = e_2 \in \bO$). This defines an orthogonal complex structure because $|Jv| = |v|$ and $J^2 = -\id$. The corresponding Kähler form $\omega_J \in \ext^2 \R^8$ is defined by
\[
\omega_J(v,w) \coloneqq \langle Jv, w\rangle = \langle v \bi, w \rangle = \langle \bi, \bar{v}w \rangle .
\]
The final equality follows from \eqref{eq:algebra_basics}. In the standard basis $e_1, \dots, e_8$ of $\bO$, the complex structure $J$ is determined by the multiplication rules in \eqref{eq:octonion_rules} and satisfies
\begin{equation*}
Je_1 = e_2, \quad Je_3 = -e_4, \quad Je_5 = -e_6, \quad Je_7 = e_8 .
\end{equation*}
Consequently, the Kähler form can be written as
\[
\omega_J = E_{12} - E_{34} - E_{56} + E_{78} .
\]
Using the complex forms
\begin{equation}
	\label{eq:complex_coordinates}
Z_1 = E_1 + \bi E_2, \quad Z_2 = E_3 - \bi E_4, \quad Z_3 = E_5 - \bi E_6, \quad Z_4 = E_7 + \bi E_8 ,
\end{equation}
the associated holomorphic volume form is given by $\Omega_J = Z_{1234}$. With this convention, the Cayley calibration can be expressed as
\begin{equation}
	\label{eq:cayley_identity}
\phi_{\mathrm{Cay}} = \Re\Omega_J - \tfrac{1}{2}\omega_J^2 .
\end{equation}
This identity follows from \cite[Proposition~1.37]{HL}. Alternatively, it can be verified directly by comparing \eqref{eq:cayley_basis} with the basis expansion of $\Re \Omega_J - \tfrac{1}{2}\omega_J^2$ as in \cite[Proposition~1.32]{HL}.

\begin{Lem}
	\label{lem:cross_i}
$\Im \Omega_J(x,y,z,w) = -\langle x \times y \times z \times w, \bi\rangle$.
\end{Lem}

\begin{proof}
Since both forms are alternating, it is enough to verify equality if $x,y,z,w$ are orthogonal. Successively applying \eqref{eq:cross}, \eqref{eq:algebra_basics}, \eqref{eq:cayley_def} and \eqref{eq:cayley_identity}, we obtain
\begin{align*}
\langle x \times y \times z \times w, \bi\rangle & = \langle \bar{x}(y(\bar{z}w)), \bi \rangle = \langle y(\bar{z}w), x\bi \rangle = \phi_{\mathrm{Cay}}(x\bi, y, z, w) \\
 & = \Re\Omega_J(Jx,y,z,w) - \tfrac{1}{2}\omega_J^2 (J x,y,z,w) \\
 & = -\Im\Omega_J(x,y,z,w) - \tfrac{1}{2}\omega_J^2 (J x,y,z,w) .
\end{align*}
Because both the cross product and $\omega_J$ are alternating,
\[
\langle x \times y \times z \times w, \bi\rangle = -\Im\Omega_J(x,y,z,w) -  \tfrac{1}{8}\mu(x,y,z,w) ,
\]
where $\mu(x,y,z,w)$ is given by
\[
\omega_J^2 (J x,y,z,w) + \omega_J^2 (x,Jy,z,w) + \omega_J^2 (x,y,Jz,w) + \omega_J^2 (x,y,z,Jw) .
\]
From
\[
\omega_J (J v, w) + \omega_J (v, J w) = - \langle v, w\rangle + \langle v, w\rangle = 0 ,
\]
it follows that $\mu = 0$, completing the proof.
\end{proof}

\subsection{Clifford algebras}
\label{subsec:clifford_algebras}

Let $(V,\langle\cdot,\cdot\rangle)$ be an oriented real inner product space of dimension $n$. The Clifford algebra $\Cl(V)$ is $\otimes V/I(V)$, where $I(V)$ is the two-sided ideal generated by the elements $x \otimes x + \langle x,x\rangle$ (see \cite[Definition~9.4]{H}). As vector spaces, $\Cl(V)$ and $\ext_\ast V$ are canonically isomorphic \cite[Proposition~9.11]{H}, with $\Cl(V)$ inheriting the inner product from $\ext_\ast V$. This inner product also provides an identification between $\Cl(V)$ and $\ext^\ast V$. The dimension of $\Cl(V)$ is $2^n$.

We will assume from now on that $n=8m$. In this case, $\Cl(V)$ and $\End(\R^{16m})$ are isomorphic as algebras \cite[Theorem~11.3]{H}. For a specific representation, the notation $\End(\bP)$ is used, where $\bP \cong \R^{16m}$ is the pinor space. The inner product $\varepsilon$ on the space of pinors is defined such that the hat involution sends every $a \in \End(\bP)$ to its adjoint, i.e., that $\varepsilon(ax,y) = \varepsilon(x,\hat ay)$ for all $x,y \in \bP$. For details on the hat involution, see \cite[Definition~9.26]{H}; and for the existence of $\varepsilon$, see \cite[Theorem~13.17]{H}. The inner product $\varepsilon$ is unique up to a conformal factor and will, in the following applications, be identified with the standard inner product with respect to a coordinate system. The space of pinors $\bP$ decomposes orthogonally into positive and negative spinors, $\bP = \bS^+ \! \oplus \bS^-$. These are the eigenspaces of the volume element $\nu \in \ext^n V$ with eigenvalues $\pm 1$. The outer product of two pinors $x,y \in \bP$ is denoted $x \scirc y \in \End(\bP) \cong \Cl(V)$, and is defined by $(x \scirc y)(z) \coloneqq \varepsilon(y,z)x$ for all $z \in \bP$.

The Clifford algebra of $\R^{n}$ with the standard inner product is denoted by $\Cl(n)$. Explicit identification for $n=8$ and $n=16$ are provided in \cite{DH} using the octonions $\bO \cong \R^8$. If $\bO^+$ and $\bO^-$ are two marked copies of $\bO$, and $\bP(8) \coloneqq \bO^+ \!\oplus \bO^-$, then $\Cl(8)$ can be identified with $\End(\bP(8))$ as an extension of
\begin{equation}
	\label{eq:clifford8model}
	\R^8 \ni v \mapsto 
	\begin{pmatrix}
		0 & R_v \\
		-R_{\bar{v}} & 0
	\end{pmatrix}
	\in \End(\bO^+ \!\oplus \bO^-) ,
\end{equation}
using the fundamental lemma of Clifford algebras \cite[Lemma~9.7]{H}. Here, $R_v \colon \bO \to \bO$, $w \mapsto wv$, denotes right multiplication with $v \in \bO$. The volume element in $\Cl(8)$ calculates, for example with Lemma~\ref{lem:octonion_product}, as
\[
	\nu \coloneqq e_1 \cdots e_8 =
	\begin{pmatrix}
		\id & 0 \\
		0 & -\id
	\end{pmatrix}
\]
and identifies the positive and negative spinors as $\bS^\pm = \bO^\pm$. For any $v \in \R^8$, represented as in \eqref{eq:clifford8model}, the following identities hold:
\begin{equation}
	\label{eq:twisted_center}
	v\nu = -\nu v \quad \text{and} \quad \nu^2 = \id .	
\end{equation}
The pinor inner product on $\bP(8) = \bS^+ \!\oplus \bS^-$ can be taken as
\[
\langle (x^+,x^-),(y^+,y^-) \rangle \coloneqq \langle x^+ , y^+ \rangle + \langle x^-, y^- \rangle ,
\]
where $\langle x, y\rangle = \Re \bar{x} y$ on the right-hand side denotes the standard inner product on $\bO \cong \R^8$ (see \cite[Lemma~14.12]{H}).

With the volume element $\nu \in \Cl(8)$ as identified above, set
\[
	E_i \coloneqq e_i \otimes \nu , \quad E_{8+i} \coloneqq \id \otimes e_i ,
\]
for $i=1, \dots, 8$. These are elements in the tensor product $\Cl(8) \otimes \Cl(8)$. Since $E_i^2 = -\id$ and $E_iE_j = - E_jE_i$ for $i \neq j$, the fundamental lemma of Clifford algebras identifies
\[
\R^{16} \cong \spa(E_1, \dots, E_{16})
\]
as a subset of
\[
\Cl(16) \cong \End(\bP(8) \otimes \bP(8)) \cong \Cl(8) \otimes \Cl(8) .
\]
Moreover, the volume element is
\[
\lambda \coloneqq \nu \otimes \nu = E_1\cdots E_{16} .
\]
Hence in this model, the positive and negative spinors in $\bP(16) \coloneqq \bP(8) \otimes \bP(8)$ are identified as
\begin{align*}
	\bS^+(16) & = (\bO^+ \!\otimes \bO^+) \oplus (\bO^- \!\otimes \bO^-) , \\
	\bS^-(16) & = (\bO^+ \!\otimes \bO^-) \oplus (\bO^- \!\otimes \bO^+) .
\end{align*}
A pinor inner product on $\bP(16)$ is induced by
\[
\langle v\otimes w, x \otimes y \rangle \coloneqq \langle v, x \rangle \langle w, y \rangle ,
\]
so that the elements $e_i^\pm \!\otimes e_j^\pm$ form an orthonormal basis. Here, $e_i^\pm$ denote the standard basis elements in $\bO^\pm$.

\section{Proof of Theorem~\ref{thm:main}}
\label{sec:main}

Working directly in coordinates and with the notation of Section~\ref{subsec:normed_algebras}, let $\C^8 = \C^4 \oplus \C^4$ with standard dual bases $Z_1,Z_2,Z_3,Z_4$ of $X_1 \coloneqq \C^4 \times \{0\}$ and $Z_5,Z_6,Z_7,Z_8$ of $X_2 \coloneqq\{0\} \times \C^4$. The holomorphic volume forms on $X_1$, $X_2$ and $\C^8$ are
\begin{equation*}
\begin{gathered}
\Omega_1 = Z_1 \wedge Z_2 \wedge Z_3 \wedge Z_4 \in \ext^4(\C^8,\C) , \\
\Omega_2 = Z_5 \wedge Z_6 \wedge Z_7 \wedge Z_8 \in \ext^4(\C^8,\C) , \\
\Omega = \Omega_1 \wedge \Omega_2 \in \ext^8(\C^8,\C) .
\end{gathered}
\end{equation*}
The Kähler forms are
\begin{equation*}
\begin{gathered}
\omega_1 = E_{1} \wedge E_2 + \cdots + E_{7} \wedge E_{8} \in \ext^2\C^8 , \\
\omega_2 = E_{9} \wedge E_{10} + \cdots + E_{15} \wedge E_{16} \in \ext^2\C^8 , \\
\omega = \omega_1 + \omega_2 \in \ext^2\C^8 .
\end{gathered}
\end{equation*}
In these coordinates, the form $\Phi$ in Theorem~\ref{thm:main} is the real part of 
\[
\alpha \coloneqq \Omega + \tfrac{1}{4!}\omega^4 + (\Omega_1 + \Omega_2)  \wedge \tfrac{1}{2}\omega^2 \in \ext^8 (\C^8,\C) .
\]
Because $\omega_i \wedge \Omega_i = 0$ for $i=1,2$, the form $\alpha$ is equal to
\begin{equation}
	\label{eq:alpha_formula}
\alpha = \Omega_1 \wedge \Omega_2 + \tfrac{1}{4!}(\omega_1 + \omega_2)^4 + \Omega_1 \wedge \tfrac{1}{2}\omega_2^2 + \tfrac{1}{2}\omega_1^2 \wedge \Omega_2  .
\end{equation}
Note that $\beta \wedge \gamma = \gamma \wedge \beta$ if $\beta$ and $\gamma$ are forms of even degree. The goal is to estimate
\[
|\alpha|(\xi) \coloneqq |\Omega(\xi)| + \tfrac{1}{4!}|\omega^4(\xi)| + |((\Omega_1 + \Omega_2)  \wedge \tfrac{1}{2}\omega^2)(\xi)|
\]
for $\xi \in \bG(\C^8,8)$. Expressing $\xi$ in normal form as in \cite[Lemma~6.13]{HL}, there are a unitary basis $e_1, \dots, e_8 \in \C^8$ and angles $0 \leq \theta_1 \leq \theta_3 \leq \theta_5 \leq \frac{\pi}{2}$, $\theta_5 \leq \theta_7 \leq \pi$, such that
\begin{align*}
\xi & = e_1 \wedge (\bi e_1 \cos \theta_1 + e_2 \sin \theta_1) \wedge e_3 \wedge (\bi e_3 \cos \theta_3 + e_4 \sin \theta_3) \; \wedge \\
& \phantom{=} \; \; e_5 \wedge (\bi e_5 \cos \theta_5 + e_6 \sin \theta_5) \wedge e_7 \wedge (\bi e_7 \cos \theta_7 + e_8 \sin \theta_7) \\
& = \xi_1 \wedges \xi_8 .
\end{align*}
Each of the three components of $|\alpha|(\xi)$ is analyzed separately. Let $0 \leq \varphi < 2\pi$ be the phase such that
\[
e^{\bi \varphi} = \det [e_1, \dots, e_8] = \Omega(e_1 \wedges e_8).
\]
First, since $e_j$ and $\bi e_j$ are linearly dependent over $\C$, we have
\begin{align}
	\nonumber
\Omega(\xi) & = \sin \theta_1 \sin \theta_3 \sin \theta_5 \sin \theta_7 \, \det [e_1, \dots, e_8] \\
	\label{eq:Omega_estimate}
 & = \sin \theta_1 \sin \theta_3 \sin \theta_5 \sin \theta_7 \, e^{\bi \varphi}.
\end{align}
Since $\omega(e_j,e_k) = \omega(e_j,\bi e_k) = 0$ for $j \neq k$, and since $\tfrac{1}{4!} \omega^4$ calibrates the complex $4$-plane
\[
\eta \coloneqq e_1 \wedge \bi e_1 \wedge e_3 \wedge \bi e_3 \wedge e_5 \wedge \bi e_5 \wedge e_7 \wedge \bi e_7 \in \bG(\C^8,8)
\]
as stated in \cite[\S1.8.2]{F}, it follows that
\begin{align}
	\nonumber
\tfrac{1}{4!}\omega^4(\xi) & = \cos \theta_1 \cos \theta_3 \cos \theta_5 \cos \theta_7 \tfrac{1}{4!} \omega^4(\eta) \\
	\label{eq:omega_estimate}
 & = \cos \theta_1 \cos \theta_3 \cos \theta_5 \cos \theta_7 .
\end{align}
We now turn to the mixed term of $|\alpha|(\xi)$. For this purpose, let $\pi_i \colon \C^8 \to X_i$ for $i=1,2$ denote the orthogonal projections, and define $f_j \coloneqq \pi_1(e_j)$ and $g_j \coloneqq \pi_2(e_j)$, so that $e_j = f_j + g_j$. By the definition of the wedge product, we have
\begin{equation*}
(\Omega_1 \wedge \omega^2)(\xi) = \sum_{I \in \Lambda(8,4)} \sgn(\sigma_I)\, \Omega_1(\xi_{I_1} \wedges \xi_{I_4})\, \omega^2(\xi_{I^c_1} \wedges \xi_{I^c_4}) .
\end{equation*}
Since $\omega(e_j,e_k) = \omega(e_j,\bi e_k) = 0$ for $j \neq k$, it follows that
\[
\omega^2(\xi_{I^c_1} \wedges \xi_{I^c_4}) = 0
\]
unless $I^c$ (or equivalently $I$) consists of pairs $(2j-1,2j)$. If $I$ is made up of such pairs, then $\sgn(\sigma_I) = 1$. Let $\Lambda$ denote the collection of all strictly increasing sequences $\lambda \colon \{1,2\} \to \{1,3,5,7\}$. For each $\lambda \in \Lambda$, let $\lambda^c \colon \{1,2\} \to \{1,3,5,7\}$ be the complementary sequence. For any $\lambda \in \Lambda$, we introduce the notation
\begin{align*}
\xi_\lambda & \coloneqq \xi_{\lambda_1} \wedge \xi_{\lambda_1 + 1} \wedge \xi_{\lambda_2} \wedge \xi_{\lambda_2 + 1} \\
 & = e_{\lambda_1} \wedge (\bi e_{\lambda_1} \cos \theta_{\lambda_1} + e_{\lambda_1+1} \sin \theta_{\lambda_1}) \wedge e_{\lambda_2} \wedge (\bi e_{\lambda_2} \cos \theta_{\lambda_2} + e_{\lambda_2+1} \sin \theta_{\lambda_2})
\end{align*}
so that
\begin{equation}
	\label{eq:mixed_eval}
(\Omega_1 \wedge \omega^2)(\xi) = \sum_{\lambda \in \Lambda} \Omega_1(\xi_{\lambda})\, \omega^2(\xi_{\lambda^c}) .
\end{equation}
For $\lambda \in \Lambda$, we have
\begin{align}
	\nonumber
\tfrac{1}{2}\omega^2(\xi_{\lambda}) & = \cos \theta_{\lambda_1} \cos \theta_{\lambda_2} \tfrac{1}{2}\omega^2(e_{\lambda_1} \wedge \bi e_{\lambda_1} \wedge e_{\lambda_2} \wedge \bi e_{\lambda_2}) \\
	\label{eq:mixed_eval2}
 & =  \cos \theta_{\lambda_1} \cos \theta_{\lambda_2} ,
\end{align}
where, in the last step, we use the fact that $e_{\lambda_1} \wedge \bi e_{\lambda_1} \wedge e_{\lambda_2} \wedge \bi e_{\lambda_2}$ spans a complex $2$-plane, which is therefore calibrated by $\frac{1}{2}\omega^2$. Since $f_j$ and $\bi f_j$ are linearly dependent over $\C$, we obtain
\begin{align}
 \nonumber
\Omega_1(\xi_{\lambda}) & = \Omega_1(\pi_{1} \xi_{\lambda_1} \wedge \pi_1 \xi_{\lambda_1+1} \wedge \pi_{1} \xi_{\lambda_2} \wedge \pi_{1} \xi_{\lambda_2+1}) \\
 \nonumber
 & = \det [\pi_{1} \xi_{\lambda_1}, \pi_1 \xi_{\lambda_1+1}, \pi_{1} \xi_{\lambda_2}, \pi_{1} \xi_{\lambda_2+1}] \\
\label{eq:mixed_eval3}
 & = \sin \theta_{\lambda_1} \sin \theta_{\lambda_2} \det(F_\lambda) ,
\end{align}
where $F_\lambda \coloneqq [f_{\lambda_1}, f_{\lambda_1+1}, f_{\lambda_2}, f_{\lambda_2+1}]$ is a complex $4 \times 4$ matrix. Similarly,
\[
\Omega_2(\xi_{\lambda}) = \sin \theta_{\lambda_1} \sin \theta_{\lambda_2}\det(G_\lambda),
\]
where $G_\lambda \coloneqq [g_{\lambda_1}, g_{\lambda_1+1}, g_{\lambda_2}, g_{\lambda_2+1}]$. 
Combining \eqref{eq:mixed_eval}, \eqref{eq:mixed_eval2} and \eqref{eq:mixed_eval3}, it follows that
\begin{align}
	\nonumber
((\Omega_1 & + \Omega_2) \wedge \tfrac{1}{2} \omega^2)(\xi) \\
	\nonumber
& = \sum_{\lambda \in \Lambda} (\Omega_1(\xi_{\lambda}) + \Omega_2(\xi_\lambda))\, \tfrac{1}{2}\omega^2(\xi_{\lambda^c}) \\
	\label{eq:mixed evaluation}
& = \sum_{\lambda \in \Lambda} \sin \theta_{\lambda_1} \sin \theta_{\lambda_2} \cos \theta_{\lambda_1^c} \cos \theta_{\lambda_2^c}( \det(F_\lambda) +  \det(G_\lambda)) .
\end{align}
The essential technical analysis of this expression is captured in the following lemma. Recall that
\[
\sigma \coloneqq Z_1 \wedge Z_2 +  Z_3 \wedge Z_4 +  Z_5 \wedge Z_6  + Z_7 \wedge Z_8 \in \ext^2 (\C^8,\C)
\]
is the standard $\C$-symplectic form on $\C^8$.

\begin{Lem}
	\label{lem:mixed_estimates}
The identity
\begin{equation}
	\label{eq:matrix_identity}
\det(G_{\lambda}) = e^{\bi\varphi} \overline{\det(F_{\lambda^c})}
\end{equation}
holds for every  $\lambda \in \Lambda$. Define
\begin{equation*}
\begin{gathered}
\Lambda' \coloneqq \{(1,3), (1,5), (1,7)\} \subset \Lambda , \\
s_j \coloneqq \sin \theta_j \quad \text{and} \quad c_j \coloneqq |\cos \theta_j| , \\
M(\theta) \coloneqq \max_{\lambda \in \Lambda'} s_{\lambda_1} s_{\lambda_2} c_{\lambda_1^c} c_{\lambda_2^c} + c_{\lambda_1} c_{\lambda_2} s_{\lambda_1^c} s_{\lambda_2^c} , \\
|\beta|(F) \coloneqq \sum_{\lambda \in \Lambda'} \bigl| \det(F_\lambda) + e^{\bi\varphi} \overline{\det(F_{\lambda^c})} \bigr| ,
\end{gathered}
\end{equation*}
where $F \coloneqq [f_1, \dots, f_8]$ is the complex $4 \times 8$ matrix composed of the first $4$ rows of $[e_1, \dots, e_8]$. Then, it follows
\begin{equation}
\label{eq:mixed_estimate}
|((\Omega_1 + \Omega_2) \wedge \tfrac{1}{2} \omega^2)(\xi)| \leq M(\theta)\, |\beta|(F) ,
\end{equation}
and
\begin{equation}
\label{eq:symplectic_estimate}
|\beta|(F) = \Re\tfrac{1}{2}\sigma^2(r_1 \wedge r_2 \wedge r_3 \wedge r_4) ,
\end{equation}
for a unitary system $r_1,r_2,r_3,r_4 \in \C^8$.
\end{Lem}

\begin{proof}
For $\lambda \in \Lambda$, define the block matrix
\[
M_\lambda \coloneqq
\begin{pmatrix}
	F_\lambda & F_{\lambda^c} \\
	G_\lambda & G_{\lambda^c}
\end{pmatrix}
.
\]
Up to a permutation with signature $1$, the columns of $M_\lambda$ are $e_1, \dots, e_8$. Therefore, it follows that $\det(M_\lambda) = e^{\bi\varphi}$, and $M_\lambda$ is unitary,
\begin{align*}
F_\lambda F_\lambda^\ast + F_{\lambda^c}F_{\lambda^c}^\ast = G_\lambda G_\lambda^\ast + G_{\lambda^c}G_{\lambda^c}^\ast = \id, \quad F_\lambda G_\lambda^\ast + F_{\lambda^c}G_{\lambda^c}^\ast = 0.
\end{align*}
Hence 
\begin{align*}
	\begin{pmatrix}
		F_\lambda & F_{\lambda^c} \\
		0 & \id
	\end{pmatrix}
	M_\lambda^\ast
=
\begin{pmatrix}
	F_\lambda & F_{\lambda^c} \\
	0 & \id
\end{pmatrix}
\begin{pmatrix}
	F_\lambda^\ast & G_\lambda^\ast \\
	F_{\lambda^c}^\ast & G_{\lambda^c}^\ast
\end{pmatrix}
=
\begin{pmatrix}
	\id & 0 \\
	F_{\lambda^c}^\ast & G_{\lambda^c}^\ast
\end{pmatrix}
\end{align*}
and
\begin{equation*}
\det(F_\lambda) e^{-\bi\varphi} =  \det(G_{\lambda^c}^\ast) =  \overline{\det(G_{\lambda^c})} .
\end{equation*}
This shows \eqref{eq:matrix_identity}. Introduce the notation
\[
z(\lambda) \coloneqq \det(F_\lambda) +  \det(G_\lambda) = \det(F_\lambda) +  e^{\bi\varphi} \overline{\det(F_{\lambda^c})}
\]
for all $\lambda \in \Lambda$. For complementary sequences $\lambda,\lambda^c \in \Lambda$, these values are related by
\begin{align*}
e^{-\bi\varphi/2} z(\lambda) & = e^{-\bi\varphi/2} \bigl( \det(F_\lambda) + e^{\bi\varphi} \overline{\det(F_{\lambda^c})} \bigr) \\
 & =
 \bigl( e^{-\bi\varphi/2} \det(F_\lambda) + \overline{e^{-\bi\varphi/2}\det(F_{\lambda^c})} \bigr) \\
 & = \overline{e^{-\bi\varphi/2} z(\lambda^c)} .
\end{align*}
In particular, $|z(\lambda)| = |z(\lambda^c)|$. Note that every element of $\Lambda$ is either in $\Lambda'$ or is the complement of an element in $\Lambda'$. Therefore, the mixed wedge product, as expressed in \eqref{eq:mixed evaluation}, is bounded by
\begin{align*}
\nonumber
|((\Omega_1 + \Omega_2) \wedge \tfrac{1}{2} \omega^2)(\xi)| & \leq \sum_{\lambda \in \Lambda'} (s_{\lambda_1} s_{\lambda_2} c_{\lambda_1^c} c_{\lambda_2^c} + c_{\lambda_1} c_{\lambda_2} s_{\lambda_1^c} s_{\lambda_2^c}) |z(\lambda)| \\
& \leq M(\theta)\, |\beta|(F) .
\end{align*}
This proves \eqref{eq:mixed_estimate}.

In order to establish \eqref{eq:symplectic_estimate}, note that $F$ is a unitary $4 \times 8$ matrix, meaning that its rows $q_1,q_2,q_3,q_4 \in \C^8$ form a unitary system. With respect to the standard dual basis of $\C^8$, we have
\[
\det(F_\lambda)
= \det\bigl([f_{\lambda_1}, f_{\lambda_1+1}, f_{\lambda_2}, f_{\lambda_2+1}]\bigr)
= (Z_{\lambda_1} \wedge Z_{\lambda_1+1} \wedge Z_{\lambda_2} \wedge Z_{\lambda_2+1})(F),
\]
where the right-hand side evaluation on $F$ is an abbreviation for the evaluation on the row-space $q_1 \wedge q_2 \wedge q_3 \wedge q_4$ of $F$. Thus explicitly,
\begin{align*}
|\beta|(F) & = \bigl|
e^{-\bi\varphi/2} Z_{1234}(F) + \overline{e^{-\bi\varphi/2} Z_{5678}(F)} \bigr| + \\
& \phantom{=} \;\; \bigl|
e^{-\bi\varphi/2} Z_{1256}(F) + \overline{e^{-\bi\varphi/2} Z_{3478}(F)} \bigr| + \\
& \phantom{=} \;\; \bigl|
e^{-\bi\varphi/2} Z_{1278}(F) + \overline{e^{-\bi\varphi/2} Z_{3456}(F)} \bigr|
\end{align*}
If $F_1$ is obtained by multiplying $F$ with $e^{-\bi\varphi/8}$, then $F_1$ is unitary and $|\beta|(F)$ is equal to
\begin{align*}
\bigl| Z_{1234}(F_1) \! + \! \overline{Z_{5678}(F_1)} \bigr| \! + \! \bigl| Z_{1256}(F_1) \! + \! \overline{Z_{3478}(F_1)} \bigr| \! + \! \bigl| Z_{1278}(F_1) \! + \! \overline{Z_{3456}(F_1)} \bigr| .
\end{align*}
If $F_2$ is obtained from $F_1$ by multiplying columns $1$ and $3$ by $e^{\bi a/2}$ and columns $5$ and $7$ by $e^{-\bi a/2}$ for $a \in \R$, then $F_2$ is unitary and
\begin{align*}
	Z_{1234}(F_2) + \overline{Z_{5678}(F_2)} & = e^{\bi a} \bigl(Z_{1234}(F_1) + \overline{Z_{5678}(F_1)} \bigr) , \\
	Z_{1256}(F_2) + \overline{Z_{3478}(F_2)} & = Z_{1256}(F_1) + \overline{Z_{3478}(F_1)} , \\
	Z_{1278}(F_2) + \overline{Z_{3456}(F_2)} & = Z_{1278}(F_1) + \overline{Z_{3456}(F_1)} .
\end{align*}
By proceeding with different columns, we can adjust the phases of all three terms independently to obtain a unitary $4 \times 8$ matrix $F'$ such that $|\beta|(F)$ is equal to
\begin{align*}
& \Re \bigl(Z_{1234}(F') \! + \! \overline{Z_{5678}(F')} \! + \!	Z_{1256}(F') \! + \! \overline{Z_{3478}(F')} \! + \!
Z_{1278}(F') \! + \! \overline{Z_{3456}(F')} \bigr) \\
& \qquad = \Re(Z_{1234} \! + \! Z_{5678} \! + \!
Z_{1256} \! + \! Z_{3478} \! + \!
Z_{1278} \! + \! Z_{3456})(F') .
\end{align*}
This is exactly $\Re\frac{1}{2}\sigma^2(r_1 \wedge r_2 \wedge r_3 \wedge r_4)$, where $r_i$ denotes the $i$th row of $F'$, thereby showing \eqref{eq:symplectic_estimate}.
\end{proof}

We now assemble the final arguments to conclude the proof of Theorem~\ref{thm:main}.

\begin{proof}[Proof of Theorem~\ref{thm:main}]
The form $\Re\tfrac{1}{2}\sigma^2 \in \ext^4 \C^8$ is a calibration, as shown in \cite[Theorem~2.38]{BH}, and $r_1 \wedge r_2 \wedge r_3 \wedge r_4$ belongs to $\bG(\C^8,4)$. Consequently, we obtain $|\beta|(F) \leq 1$. Notably, among the values $\sin \theta_i$ and $\cos \theta_i$, only $\cos \theta_7$ can be negative. Let $(i,j,k)$ be a permutation of $(1,3,5)$ such that
\[
M(\theta) = s_{i} s_{j} c_{k} c_{7} + c_{i} c_{j} s_{k} s_{7} .
\]
By combining \eqref{eq:Omega_estimate}, \eqref{eq:omega_estimate} and \eqref{eq:mixed_estimate}, we obtain
\begin{align*}
|\alpha|(\xi) & \leq s_i s_j s_k s_7 + c_i c_j c_k c_7 + s_{i} s_{j} c_{k} c_{7} + c_{i} c_{j} s_{k} s_{7} \\
 & = (c_k c_7 + s_k s_7)(c_i c_j + s_i s_j) \\
 & = \delta \cos(\theta_7 - \delta\theta_k) \cos(\theta_j - \theta_i) ,
\end{align*}
where $\delta = 1$ if $\theta_7 \leq \frac{\pi}{2}$ and $\delta = -1$ otherwise. In any case, we have $|\alpha|(\xi) \leq 1$ with equality occurring only if the angles form pairs. That is, there exists a partition $\{a,b\} \cup \{c,d\} = \{1,3,5,7\}$ such that
\begin{equation}
\label{eq:angle_condition}
\theta_a = \theta_b \quad \text{and} \quad \sin \theta_c = \sin \theta_d .
\end{equation}
The second condition is equivalent to $\theta_c = \theta_d$ or $\theta_c + \theta_d = \pi$.

We now return to the setting of Theorem~\ref{thm:main}. By making suitable coordinate choices to identify $X_1$ and $X_2$ with $\C^4$, it follows from the bound $|\alpha|(\xi) \leq 1$ just established that we have the estimate
\[
|\Omega(\xi)| + |\tfrac{1}{4!}\omega^4(\xi)| + |(\Omega'_1 \wedge \tfrac{1}{2}\omega^2)(\xi) + (\Omega'_2 \wedge \tfrac{1}{2}\omega^2)(\xi)| \leq 1
\]
for all $\xi \in \bG(X,8)$ and for arbitrary holomorphic volume forms $\Omega'_1$ on $X_1$ and $\Omega'_2$ on $X_2$. Note that $|\Omega(\xi)| = |(\Omega_1' \wedge \Omega_2')(\xi)|$ for any holomorphic volume form on $X = X_1 \oplus X_2$. Let $\Omega_1$ and $\Omega_2$ be the fixed holomorphic volume forms in Theorem~\ref{thm:main}. Since $e^{\bi t}\Omega_1$ is again a holomorphic volume form for any $t \in \R$, we can apply the estimate to $\Omega'_1 = e^{\bi t}\Omega_1$ and $\Omega'_2 = \Omega_2$. We may choose $t$ (depending on $\xi$) to align the phases such that
\[
|(e^{\bi t}\Omega_1 \wedge \tfrac{1}{2}\omega^2)(\xi) + (\Omega_2 \wedge \tfrac{1}{2}\omega^2)(\xi)\bigr|
= |(\Omega_1 \wedge \tfrac{1}{2}\omega^2)(\xi)| + |(\Omega_2 \wedge \tfrac{1}{2}\omega^2)(\xi)| .
\]
This proves the inequality in Theorem~\ref{thm:main} with the four absolute values.

For the remainder of the proof, it is useful to express $\Phi$ as in \eqref{eq:alpha_formula},
\begin{equation}
	\label{eq:phi_count}
\Phi = \Re \Omega + \tfrac{1}{4!}\omega^4 + \Re\Omega_1 \wedge \tfrac{1}{2}\omega_2^2 + \tfrac{1}{2}\omega_1^2 \wedge \Re\Omega_2 .
\end{equation}
By Lemma~\ref{lem:standard_calibrations}, the forms $\Re\Omega$ and $\tfrac{1}{4!}\omega^4$ are self-dual. Let $\ast_i$ denote the Hodge star operator on $\ext^\ast X_i$ for $i = 1,2$. Since $X_1$ and $X_2$ are orthogonal, we have
\begin{align*}
\ast\bigl(\Re\Omega_1 \wedge \tfrac{1}{2}\omega_2^2\bigr) & = \bigl(\ast_1\!\Re\Omega_1\bigr) \wedge \bigl(\ast_2 \tfrac{1}{2} \omega_2^2\bigr) = \Re\Omega_1 \wedge \tfrac{1}{2}\omega_2^2
\end{align*}
and the same holds for $\tfrac{1}{2}\omega_1^2 \wedge \Re\Omega_2$. Therefore, $\Phi$ is self-dual.

Let $\Lambda \subset \Lambda(16,8)$ be the set of sequences $I$ for which $E_I \in \ext^8 \R^{16}$ contributes nontrivially to $\Phi$. By Lemma~\ref{lem:standard_calibrations}, we can express
\[
\Phi = \sum_{I \in \Lambda} \lambda_I E_I
\]
with $\lambda_I = \pm 1$. Let $\Lambda_1, \Lambda_2, \Lambda_3, \Lambda_4 \subset \Lambda$ correspond to the contributions from the four components of $\Phi$ in \eqref{eq:phi_count}, listed in order from left to right. These sets are pairwise disjoint, as follows from their combinatorial descriptions, derived similarly to the proof of Lemma~\ref{lem:standard_calibrations}. This disjointness, also follows from the comass bound which has already been established. By Lemma~\ref{lem:standard_calibrations}, the cardinalities of these sets are given by
\begin{equation*}
|\Lambda_1| = 2^7  = 128, \quad |\Lambda_2| = \binom{8}{4} = 70, \quad |\Lambda_3| = |\Lambda_4| = 2^3 \cdot \binom{4}{2} = 48.
\end{equation*}
Thus the total count is $|\Lambda| = 128 + 70 + 96 = 294$. Because $\Phi$ is self-dual,
\[
294 \vol = \langle \Phi, \Phi \rangle \vol = \Phi \wedge \! \ast\Phi = \Phi^2 .
\]
This completes the proof of Theorem~\ref{thm:main}.
\end{proof}

The quantity
\[
\mathrm{Wirt}_n \coloneqq \sup \left\{\frac{|\alpha^2|}{\|\alpha\|^2} : \alpha \in \ext^n \R^{2n}, \alpha \neq 0\right\}
\]
was introduced by Bangert, Katz, Shnider and Weinberger in \cite[Definition~5.1]{BKSW}. More recently, Hebda and Katz defined
\[
C_{n,k} \coloneqq \sup \left\{\frac{|\alpha|}{\|\alpha\|} : \alpha \in \ext^k \R^n, \alpha \neq 0\right\}
\]
in \cite[Definition~2.1]{HK} and refined the upper bound $C_{n,k}^2 \leq \binom{n}{k}$ from \cite[\S1.8.1]{F}. These quantities serve as upper bounds for stable systolic ratios, as stated in \cite[Proposition~5.4]{BKSW}, \cite[Theorem~5.1]{HK} and \cite{K}. Using the same argument as in the proof of  \cite[Lemma~5.3]{BKSW}, it is clear that $\mathrm{Wirt}_n \leq C_{2n,n}^2$:
\[
|\alpha^2| = |\langle \alpha, \! \ast \alpha \rangle| \leq |\alpha||\!\ast\! \alpha| = |\alpha|^2.
\]
Except for the trivial cases $C_{n,1} = C_{1,n} = 1$, known values are
\begin{itemize}
	\item $\mathrm{Wirt}_2 = C_{4,2}^2 = 2$, achieved by Kähler forms, and more generally, $C_{n,2}^2 = C_{n,n-2}^2 = \lfloor\frac{n}{2}\rfloor$ (see  \cite[Remark~2.2]{HK});
	\item $\mathrm{Wirt}_3 = C_{6,3}^2 = 4$ (see \cite[Proposition~3.1]{HK}); this builds on the classification of calibrations in $\ext^3 \R^6$ by Morgan \cite{M};
	\item $\mathrm{Wirt}_4 = C_{8,4}^2 = 14$, achieved by Cayley calibrations (see \cite[Proposition~9.1]{BKSW});
	\item $C_{7,3}^2 = C_{7,4}^2 = 7$ see (\cite[Proposition 3.3]{HK}).
\end{itemize}
Note that all the values in middle dimension listed above do not satisfy the strict inequality \eqref{eq:main_estimate}. Theorem~\ref{thm:main} provides the lower bound $294 \leq \mathrm{Wirt}_8 \leq C_{16,8}^2$.

\section{Calibrated planes}

The characteristic angles of calibrated planes in Theorem~\ref{thm:main2} exhibit symmetries similar to those arising in the case of Cayley calibrations. The following proposition is essentially adapted from \cite[Proposition~1.34]{HL}.

\begin{Prop}[\cite{HL}]
	\label{prop:cayley_normal}
	Consider the $4$-form on $\C^4$ defined by
	\[
	\phi_{\rm{C}} \coloneqq \Re\Omega + \tfrac{1}{2}\omega^2 \in \ext^4\C^4 \cong \ext^4 \R^8 .
	\]
	The form $\phi_{\rm{C}}$ is a self-dual Cayley calibration. Moreover, the inequality
	\[
	|\Omega(\xi)| + |\tfrac{1}{2}\omega^2(\xi)| \leq 1
	\]
	holds for all $\xi \in \bG(\C^4,4)$, and $\xi$ is calibrated by $\phi_{\rm{C}}$ if and only if it can be written as
	\begin{equation}
		\label{eq:cayley_planes}
		\xi = e_1 \wedge (\bi e_1 \cos \theta + e_2 \sin \theta) \wedge e_3 \wedge (\bi e_3 \cos \theta + e_4 \sin \theta)
	\end{equation}
	for a special unitary basis $e_1,e_2,e_3,e_4 \in \C^4$ (i.e., $[e_1,e_2,e_3,e_4] \in \operatorname{SU}(4)$) and $0 \leq \theta \leq \frac{\pi}{2}$.
\end{Prop}

\begin{proof}[Sketch of a proof]
Write $\xi \in \bG(\C^4,4)$ in general normal form with respect to a unitary basis $e_1,e_2,e_3,e_4 \in \C^4$ and angles $0 \leq \theta_1 \leq \frac{\pi}{2}$, $\theta_1 \leq \theta_2 \leq \pi$ as in \cite[Lemma~6.13]{HL}. Then
\[
|\Omega(\xi)| + |\tfrac{1}{2}\omega^2(\xi)| = \sin \theta_1 \sin \theta_2 + \cos \theta_1 |\cos \theta_2| \leq 1,
\]
with equality if and only if $\theta_1 = \theta_2$ or $\theta_1 + \theta_2 = \pi$. If $\phi_{\mathrm{C}}(\xi) = 1$, then $\theta_1 = \theta_2$. In case $\theta_1 = 0$, then up to changing the phase of $e_2$, which does not affect $\xi$, we can assume $\det [e_1,e_2,e_3,e_4] = 1$. This condition is automatically satisfied when $\theta_1 > 0$.
\end{proof}

The goal is to establish Theorem~\ref{thm:main2}. We continue using the notation from the previous section. Assume that $\xi \in \bG(\C^8,8)$ is calibrated by $\Phi$, meaning $\Phi(\xi) = 1$, and is given in normal form with angles $\theta_1,\theta_3,\theta_5,\theta_7$ as before. Since $\Phi(\xi) = 1$ implies $|\alpha|(\xi) = 1$, the angles must form pairs as in \eqref{eq:angle_condition}. Furthermore, 
\begin{align*}
\Phi(\xi) & = \sin \theta_1 \sin \theta_3 \sin \theta_5 \sin \theta_7 \, \Re e^{\bi \varphi} + \cos \theta_1 \cos \theta_3 \cos \theta_5 \cos \theta_7 \\
& \phantom{=} \; + \sum_{\lambda \in \Lambda} \sin \theta_{\lambda_1} \sin \theta_{\lambda_2} \cos \theta_{\lambda_1^c} \cos \theta_{\lambda_2^c} \Re\, ( \det(F_\lambda) + \det(G_{\lambda}) ) .
\end{align*}
To classify the calibrated planes, we consider different cases guided by the following two lemmas.

\begin{Lem}
	\label{lem:split_case}
If
\[
|\det(F_\lambda) + \det(G_\lambda)| = 1
\]
for some $\lambda \in \Lambda$, then $\xi = \xi_1 \wedge \xi_2$ with $\xi_i \in \bG(X_i,4)$. In particular, this condition holds if $|\alpha|(\xi) = 1$ and there exist indices $i$ and $j$ such that $\sin(\theta_i) \neq \sin(\theta_j)$ and $\sin(\theta_i) \in \{0,1\}$.
\end{Lem}

\begin{proof}
Without loss of generality we assume that $\lambda = (1,3)$. Using Hadamard's inequality and the Cauchy–Schwarz inequality, we obtain
\begin{align}
\nonumber
1 & = |\det(F_\lambda) +  \det(G_\lambda)| \\
\nonumber
& \leq |\det [f_1,f_2,f_3,f_4]| + |\det [g_1,g_2,g_3,g_4]| \\
\nonumber
& \leq |f_1||f_2||f_3||f_4| + |g_1||g_2||g_3||g_4| \\
\nonumber
& \leq \frac{1}{6} \sum_{1 \leq i < j \leq 4} |f_i||f_j| + |g_i||g_j| \\
\nonumber
& \leq \frac{1}{6} \sum_{1 \leq i < j \leq 4} (|f_i|^2 + |g_i|^2)^\frac{1}{2} (|f_j|^2 + |g_j|^2)^\frac{1}{2} \\
\label{eq:extreme_case_estimate}
& = 1 .
\end{align}
Equalities hold only if $|f_i| = 1$ and $|g_i| = 0$ for all $i$, or $|f_i| = 0$ and $|g_i| = 1$ for all $i$. Consequently, $e_1, e_2, e_3, e_4$ form a unitary basis of either $X_1$ or $X_2$. Without loss of generality, we assume the former. Then $e_5, e_6, e_7, e_8$ form a unitary basis of $X_2$, so $\xi$ decomposes as claimed.

From \eqref{eq:angle_condition}, we know that $\sin(\theta_i)$, for $i = 1, 2, 3, 4$, take at most two distinct values. Now, suppose that there exist indices $i$ and $j$ such that $\sin(\theta_i) \neq \sin(\theta_j)$ and $\sin(\theta_i) \in \{0,1\}$. In this case, there is exactly one $\lambda \in \Lambda$ for which
\[
s_{\lambda_1}s_{\lambda_2} c_{\lambda_1^c} c_{\lambda_2^c} \neq 0 .
\]
For $|\alpha|(\xi) = 1$ to hold, equality must be attained in
\[
|((\Omega_1 + \Omega_2) \wedge \tfrac{1}{2} \omega^2)(\xi)| \leq s_{\lambda_1}s_{\lambda_2} c_{\lambda_1^c} c_{\lambda_2^c}.
\]
This follows from \eqref{eq:mixed_estimate} and imposes the condition $|\det(F_\lambda) + \det(G_\lambda)| = 1$, completing the proof of the lemma.
\end{proof}

For each $i=1,2$, define the form
\[
\alpha_i \coloneqq \Omega_i + \tfrac{1}{2}\omega_i^2 \in \ext^4(X_i,\C) \subset \ext^4(\C^8,\C),
\]
the real part of which is the Cayley calibration as described in Proposition~\ref{prop:cayley_normal}.

\begin{Lem}
	\label{lem:calibrated_split_planes}
If $\xi = \xi_1 \wedge \xi_2$ with $\xi_i \in \ext_4 X_i$, then
\[
\alpha(\xi) = \alpha_1(\xi_1) \alpha_2(\xi_2) .
\]
If $\xi = \xi_1 \wedge \xi_2$ for $\xi_i \in \bG(X_i,4)$ and $\Phi(\xi) = 1$, then one of the following cases applies:
\begin{enumerate}
	\item $\xi$ is a special Lagrangian $8$-plane.
	\item Either $\xi_1$ and $\xi_2$, or $-\xi_1$ and $-\xi_2$, are Cayley $4$-planes as in \eqref{eq:cayley_planes}.
\end{enumerate}
\end{Lem}

\begin{proof}
From \eqref{eq:alpha_formula}, we obtain
\[
\alpha = \Omega_1 \wedge \Omega_2 + \Omega_1 \wedge \tfrac{1}{2}\omega_2^2 + \tfrac{1}{2}\omega_1^2 \wedge \Omega_2 + \frac{1}{4!}\sum_{k = 0}^4 \binom{4}{k} \omega_1^k \wedge \omega_2^{4-k} .
\]
For the first three terms, the evaluation on $\xi_1 \wedge \xi_2$ is straightforward. In the remaining sum, only $\omega_1^2 \wedge \omega_2^2$ contributes with
\[
\tfrac{6}{4!}(\omega_1^2 \wedge \omega_2^2)(\xi_1 \wedge \xi_2) = \tfrac{1}{2} \omega_1^2(\xi_1) \tfrac{1}{2}\omega_2^2(\xi_2) .
\]
This establishes that $\alpha(\xi_1 \wedge \xi_2) = \alpha_1(\xi_1)\alpha_2(\xi_2)$.

Let $\xi_i \in \bG(X_i,4)$ for $i=1,2$ and assume that $\Phi(\xi) = 1$. Then, $|\alpha_1(\xi_1)| = |\alpha_2(\xi_2)| = 1$ by the first part. Expressed in normal form,
\begin{align*}
\xi_1 & = f_1 \wedge (\bi f_1 \cos \theta_1 + f_2 \sin \theta_1) \wedge f_3 \wedge (\bi f_3 \cos \theta_2 + f_4 \sin \theta_2) , \\
\xi_2 & = g_1 \wedge (\bi g_1 \cos \theta_3 + g_2 \sin \theta_3) \wedge g_3 \wedge (\bi g_3 \cos \theta_4 + g_4 \sin \theta_4) ,
\end{align*}
with respect to unitary bases $f_1,f_2,f_3,f_4 \in X_1$ and $g_1,g_2,g_3,g_4 \in X_2$, where the characteristic angles satisfy $0 \leq \theta_1,\theta_3 \leq \frac{\pi}{2}$, $0 \leq \theta_2,\theta_4 \leq  \pi$. Define $0 \leq \psi_1,\psi_2 < 2\pi$ such that
\[
e^{\bi \psi_1} = \det [g_1, \dots, g_4] \quad \text{and} \quad e^{\bi \psi_2} = \det [f_1, \dots, f_4] .
\]
Equality $|\alpha_1(\xi_1)| = 1$ is satisfied only if $\theta_1 = \theta_2$ or $\theta_1 + \theta_2 = \pi$. Similarly, as in the proof of Proposition~\ref{prop:cayley_normal}, $\theta_3 = \theta_4$ or $\theta_3 + \theta_4 = \pi$ . Next, we classify the calibrated planes for which
\begin{align*}
1 & = \alpha(\xi) = \alpha_1(\xi_1)\alpha_2(\xi_2) \\
& = (\sin \theta_1 \sin \theta_2 e^{\bi \psi_1} + \cos \theta_1 \cos \theta_2)(\sin \theta_3 \sin \theta_4 e^{\bi \psi_2} + \cos \theta_3 \cos \theta_4) .
\end{align*}

\emph{Case A: $\psi_1 \notin \{0,\pi\}$ or $\psi_2 \notin \{0,\pi\}$.} Without loss of generality, assume $\psi_1 \notin \{0,\pi\}$. Then we have
\[
|\sin \theta_1 \sin \theta_2 e^{\bi \psi_1} + \cos \theta_1 \cos \theta_2| < |\sin \theta_1 \sin \theta_2| + |\cos \theta_1 \cos \theta_2| \leq 1 ,
\]
unless $\cos \theta_1 \cos \theta_2 = 0$. This occurs only when $\theta_1 = \theta_2 = \frac{\pi}{2}$, leading to
\[
1 = e^{\bi \psi_1}(\sin \theta_3 \sin \theta_4 e^{\bi \psi_2} + \cos \theta_3 \cos \theta_4) .
\]
Therefore, we must have $\psi_2 \notin \{0,\pi\}$, and by applying the same reasoning, it follows that $\theta_3 = \theta_4 = \frac{\pi}{2}$ and $e^{\bi(\psi_1 + \psi_2)} = 1$. Consequently, $\det{[f_1, \dots, f_4,g_1, \dots, g_4]} = 1$, meaning that $\xi$ is a special Lagrangian $8$-plane.

For the remaining cases, we assume $\psi_1,\psi_2 \in \{0,\pi\}$. This in particular implies that $\alpha_i(\xi_i)$ is $\pm 1$.

\emph{Case B: $\alpha_1(\xi_1) = 1$.} This implies $\alpha_2(\xi_2) = 1$. Thus $\xi_1$ and $\xi_2$ are Cayley $4$-planes as in Proposition~\ref{prop:cayley_normal}.

\emph{Case C: $\alpha_1(\xi_1) = -1$.} This implies also $\alpha_2(\xi_2) = -1$. Then $\xi = (-\xi_1) \wedge (-\xi_2)$ and case B applies.
\end{proof}

\begin{proof}[Proof of Theorem~\ref{thm:main2}]
\emph{Case 1: $\sin(\theta_i) = 1$ for all $i$.} In this case, $1 = \Phi(\xi) = \Re e^{\bi \varphi}$, which implies $\varphi = 0$. Thus, $\xi$ is a special Lagrangian $8$-plane.

\emph{Case 2: $\sin(\theta_i) = 0$ for all $i$.} If $\theta_i = 0$ for all $i$, then $\xi$ is a complex $4$-plane, and we obtain $1 = \Phi(\xi) = \frac{1}{4!}\omega^4(\xi)$. If $\theta_1 = \theta_3 = \theta_5 = 0$ and $\theta_7 = \pi$, then $\Phi(\xi) = -1$, meaning $\xi$ is not calibrated by $\Phi$.

\emph{Case 3: There are two indices $i,j$ with $\sin(\theta_i) \neq \sin(\theta_j)$ and $\sin(\theta_i) \in \{0,1\}$.} This case is addressed by Lemma~\ref{lem:split_case} and Lemma~\ref{lem:calibrated_split_planes}. It thus follows that $\xi$ is a special Lagrangian $8$-plane or a product of two Cayley $4$-planes.

\emph{Case 4: $\sin(\theta_i) \notin \{0,1\}$ for all $i$.} First, the product of cosines resulting from $\omega^4(\xi)$ is nonzero. Since this term contributes to $\Phi(\xi)$, the contribution has to be positive. This implies that $\theta_7 \leq \frac{\pi}{2}$. From \eqref{eq:angle_condition} we deduce that the angles form two equal pairs. Given their ordering by index, we obtain
\[
0 < \theta_1 = \theta_3 \leq \theta_5 = \theta_7 < \tfrac{\pi}{2} .
\]
In the previous section, every instance of $c_7 = |\cos \theta_7|$ can be replaced with $\cos \theta_7$. The maximal value $M(\theta)$ is attained for $\lambda = (1,3)$ and given by
\[
M(\theta) = s_1^2c_5^2 + c_1^2 s_5^2 .
\]
The corresponding trigonometric terms for $(1,5)$ and $(1,7)$ are given by $2s_1c_1s_5c_5$. If there exists some $\lambda \in \Lambda$ such that $|\det(F_\lambda) + \det(G_\lambda)| = 1$, then Lemma~\ref{lem:split_case} and Lemma~\ref{lem:calibrated_split_planes} provide a characterization of the calibrated planes. Consequently, we can assume that $|\det(F_\lambda) + \det(G_\lambda)| < 1$ for all $\lambda \in \Lambda$. Since $|\alpha|(\xi) = 1$ requires $|\beta|(F) = 1$, where $|\beta|(F)$ is defined in Lemma~\ref{lem:mixed_estimates}, there must be $\lambda \in \{(1,5),(1,7)\}$ with $|\det(F_\lambda) + \det(G_\lambda)| > 0$. With equality in \eqref{eq:mixed_estimate}, we obtain
\[
s_{1}^2 c_{5}^2 + c_{1}^2 s_{5}^2 = M(\theta) = 2s_{1} s_{5} c_{1} c_{5} .
\]
This implies that $0 = (s_1c_5 - c_1s_5)^2$, leading to $\theta_1 = \theta_5$. Consequently, all angles must be equal, say $0 < \theta < \frac{\pi}{2}$. Since $\Re e^{\bi \varphi} < 1$ for $0 < \varphi < 2\pi$, it follows that $\varphi = 0$ is a necessary condition for $\Phi(\xi) = 1$. Therefore,
\begin{align*}
1 = \Phi(\xi) & = \sin(\theta)^4 + \cos(\theta)^4 + \sin(\theta)^2\cos(\theta)^2 \sum_{\lambda \in \Lambda} \Re \bigl( \det(F_\lambda) + \overline{\det(F_{\lambda^c})} \bigr) \\
 & = \sin(\theta)^4 + \cos(\theta)^4 + 2 \sin(\theta)^2\cos(\theta)^2 \Re \sum_{\lambda \in \Lambda'} \det(F_\lambda) + \det(F_{\lambda^c}) .
\end{align*}
Since $0 < \theta < \frac{\pi}{2}$, equality is achieved only if
\[
1 = \Re \sum_{\lambda \in \Lambda'} \det(F_\lambda) + \det(F_{\lambda^c}) = \Re \sum_{\lambda \in \Lambda'} \det(G_\lambda) + \det(G_{\lambda^c}) .
\]
Let $r_i$ denote the $i$th row of the matrix $[e_1, \dots, e_8]$. Thus, both
\[
\rho_1 \coloneqq r_1 \wedge r_2 \wedge r_3 \wedge r_4 \in \bG(\C^8,4) \quad \text{and} \quad \rho_2 \coloneqq r_5 \wedge r_6 \wedge r_7 \wedge r_8 \in \bG(\C^8,4)
\]
are calibrated by
\begin{equation*}
	\label{eq:symplectic_basis}
\Re \tfrac{1}{2} \sigma^2 = \Re(Z_{1234} \! + \! Z_{5678} \! + \!
Z_{1256} \! + \! Z_{3478} \! + \!
Z_{1278} \! + \! Z_{3456}) \in \ext^4 \C^8 .
\end{equation*}
It follows from \cite[Theorem~2.38]{BH} that
\[
H_1 \coloneqq \spa_\C(r_1,r_2,r_3,r_4) \quad \text{and} \quad H_2 \coloneqq \spa_\C(r_5,r_6,r_7,r_8)
\]
are quaternionic $2$-planes, where the identification of $\C^8$ with $\bH^4$ is specified in Section~\ref{subsec:normed_algebras}. On each quaternionic $2$-plane of $\bH^4$, the form $\Re\tfrac{1}{2}\sigma^2$ is a special Lagrangian calibration. Thus, $\rho_1$ and $\rho_2$ are special Lagrangian (real) $4$-planes each lying within one of two orthogonal quaternionic $2$-planes, denoted by $H_1$ and $H_2$. 

The subspaces $X_1 = \C^4 \times \{0\}$ and $X_2 = \{0\} \times \C^4$ of $\C^8$ form a pair of orthogonal quaternionic $2$-planes. Let $S \in \mathrm{Sp}(4)$ be such that $S(X_1) = H_1$ and $S(X_2) = H_2$, and define $q_i \coloneqq S^{-1}(r_i)$ for $i=1, \dots, 8$. Since $S$ preserves $\sigma$ (and thus $\sigma^2$), it follows that
\begin{align*}
1 & = \Re\tfrac{1}{2}\sigma^2 (q_1 \wedge q_2 \wedge q_3 \wedge q_4) = \Re Z_{1234}(q_1 \wedge q_2 \wedge q_3 \wedge q_4) , \\
1 & = \Re\tfrac{1}{2}\sigma^2 (q_5 \wedge q_6 \wedge q_7 \wedge q_8) = \Re Z_{5678}(q_5 \wedge q_6 \wedge q_7 \wedge q_8) .
\end{align*}
This implies that $U_1,U_2 \in \mathrm{SU}(4)$ for $U_1 \coloneqq [q_1,q_2,q_3,q_4]^T$ and $U_2 \coloneqq [q_4,q_5,q_6,q_7]^T$. We see that the $i$th row of $U \coloneqq \operatorname{diag}(U_1,U_2)$ is $q_i$ and the $i$th row of $US^T$ is $S(q_i) = r_i$. Hence,
\[
[e_1, \dots, e_n] = U S^T .
\]
Since $S^T$ belongs to $\mathrm{Sp}(4)$, this confirms that if a calibrated plane $\xi$ with respect to $\Phi$ does not fall under cases (1), (2) or (3) of Theorem~\ref{thm:main2}, then it must belong to Case (4).

The converse direction, namely that any $\xi$ satisfying one of the four conditions is indeed calibrated by $\Phi$, is evident for (1) and (2), while (3) follows directly form Lemma~\ref{lem:calibrated_split_planes}. For (4), noting that both $U$ and $S$ lie in $\mathrm{SU}(8)$, it follows that
\[
e^{\bi \varphi} = \det [e_1, \dots, e_8] = 1 .
\]
Using the same notation as above, we obtain
\begin{align*}
\Phi(\xi) & = \sin(\theta)^4 + \cos(\theta)^4 + 2 \sin(\theta)^2\cos(\theta)^2 \Re \sum_{\lambda \in \Lambda'} \det(F_\lambda) + \det(F_{\lambda^c}) \\
 & = \sin(\theta)^4 + \cos(\theta)^4 + 2 \sin(\theta)^2\cos(\theta)^2 \Re\tfrac{1}{2}\sigma^2 (\rho) \\
 & = \sin(\theta)^4 + \cos(\theta)^4 + 2 \sin(\theta)^2\cos(\theta)^2 \\
 & = 1 ,
\end{align*}
provided that we establish $\Re\tfrac{1}{2}\sigma^2 (\rho) = 1$ for $\rho \coloneqq r_1 \wedge r_2 \wedge r_3 \wedge r_4$ spanned by the first four rows of $[e_1, \dots, e_8] = US^T$. If $q_1, \dots, q_8$ are the rows of $U$, then
\[
\Re\tfrac{1}{2}\sigma^2(q_1 \wedge q_2 \wedge q_3 \wedge q_4) = \Re(Z_{1234})(q_1 \wedge q_2 \wedge q_3 \wedge q_4) = \Re\det(U_1) = 1 .
\]
The rows are related by $r_i = S(q_i)$. Since $S$ preserves $\sigma$, we have $1 = \Re\tfrac{1}{2}\sigma^2(\rho)$. This proves Theorem~\ref{thm:main2}.
\end{proof}

We briefly describe the normal forms of (1), (2) and (3) that appear in Theorem~\ref{thm:main2}. Case (1) is also characterized by $\theta_i = \frac{\pi}{2}$ for all $i$, and 
\[
\Omega(\xi) = \det [e_1, \dots, e_8] = 1 \ .
\]
Case (2) is characterized by $\frac{1}{4!}\omega^4(\xi) = 1$, which holds if and only if $\xi$ has a normal form representation with $\theta_i = 0$ for all $i$. Case (3) means that $\xi_i \in \bG(X_i,4)$ is calibrated by $\Re \Omega_i + \tfrac{1}{2}\omega_i^2$ for $i = 1,2$ as in Proposition~\ref{prop:cayley_normal}. This occurs if and only if $\xi$ has a normal form representation for basis vectors $e_1,e_2,e_3,e_4 \in X_1$ and $e_5,e_6,e_7,e_8 \in X_2$ with
\[
\det [e_1,e_2,e_3,e_4]  = \det [e_5,e_6,e_7,e_8] = 1,
\]
and the characteristic angles satisfy $\theta_1 = \theta_2$ and $\theta_3 = \theta_4$.

\section{Spinor product realization}

Building on the preparations laid out in Section~\ref{subsec:clifford_algebras}, we now proceed to prove Theorem~\ref{thm:spinorproduct}. Define the elements
\[
s \coloneqq \bone^+ \! \otimes \bone^+, \quad s' \coloneqq \bi^+ \! \otimes \bi^+, \quad t \coloneqq s + s'
\]
of $\bO^+ \! \otimes \bO^+ \subset \bP(16)$. Observe that $\bone = e_1$ and $\bi = e_2$ are standard basis elements of $\bO \cong \R^8$,  in accordance with the conventions established in Section~\ref{subsec:normed_algebras}. In the discussion that follows, our primary objective is to calculate the homogeneous parts $\phi_4$ and $\phi_8$ of
\[
\phi \coloneqq 16^2 \, t \scirc s \in \ext^* \R^{16} .
\]
The identity
\begin{equation}
	\label{eq:phi_2n}
\phi_{2n}(X_1 \wedges X_{2n}) = \langle X_1 \cdots X_{2n} s, t\rangle
\end{equation}
holds for orthogonal vectors $X_1, \dots, X_{2n} \in \R^{16} \subset \End(\bP(16))$, as a consequence of Theorem~\ref{thm:spinorprod_dh}. We are therefore interested in evaluating the right-hand side in \eqref{eq:phi_2n}. Each $X \in \R^{16} \subset \End(\bP(16))$ admits the decomposition
\[
X = x \otimes \nu + \id \otimes y,
\]
where $x,y \in \R^8 \subset \End(\bP(8))$.

\begin{Lem}
	\label{lem:clifford_prod}
If $X_1, \dots, X_{2n} \in \R^{16}$ and $a = a_1 \otimes a_2$, $b = b_1 \otimes b_2$ for $a_1, a_2, b_1, b_2 \in \bO^+$, then
\begin{align*}
\langle X_1 \cdots X_{2n} a, b\rangle & = \sum_{k=0}^n \sum_{I \in \Lambda(2n,2k)} \sgn(\sigma_I) \bigl\langle x_{I_1} \cdots x_{I_{2k}} a_1, b_1 \bigr\rangle \bigl\langle y_{I_1^c} \cdots y_{I_{2(n-k)}^c} a_2, b_2 \bigr\rangle .
\end{align*}
Note that the products are taken in $\End(\bP(16))$ and $\End(\bP(8))$, respectively.
\end{Lem}

\begin{proof}
The product $X_1\cdots X_{2n}$ expands into a sum of terms of the form $x \otimes y$, where $x,y \in \End(\bP(8))$, and each such term evaluates to
\[
\langle (x \otimes y) a, b \rangle = \langle xa_1, b_1 \rangle \langle ya_2, b_2 \rangle .
\]
Moreover, if $x$ is formed by multiplying an odd number of $x_i$'s, then it is of anti-block-diagonal structure, as described in \eqref{eq:clifford8model}. Consequently, it maps $\bO^+$ into $\bO^-$ and the inner product $\langle xa_1, b_1 \rangle$ vanishes. For any $I \in \Lambda(2n,2k)$, define $\delta_{i} \coloneqq \nu$ if $i \in I$ and $\delta_i \coloneqq y_i$ if $i \notin I$. Since $\nu y_i = -y_i\nu$ and $\nu^2 = \id$ by \eqref{eq:twisted_center}, it follows that
\begin{align*}
\delta_1 \cdots \delta_{2n} = \sgn(\sigma_I) \nu^{2k} y_{I^c_1} \cdots y_{I^c_{2(n-k)}} = \sgn(\sigma_I) y_{I^c_1} \cdots y_{I^c_{2(n-k)}} .
\end{align*}
This proves the lemma.
\end{proof}

To express Lemma~\ref{lem:clifford_prod} in terms of octonionic products, we define the multilinear maps $P_{2k} \colon \bO^{2k} \to \bO$ recursively by $P_0 \coloneqq \bone$ and
\begin{align*}
P_{2(k+1)}(x_1, \dots, x_{2(k+1)}) \coloneqq (P_{2k}(x_{3}, \dots, x_{2(k+1)})\bar{x}_2)x_1 .
\end{align*}

\begin{Lem}
\label{lem:octonion_prod}
The following identity holds for even $n$:
\begin{equation}
	\label{eq:phi_identity}
\phi_{2n} = \sum_{k=0}^n \Alt(p_{2k,\bone}) \wedge \Alt(q_{2(n-k),\bone}) + \Alt(p_{2k,\bi}) \wedge \Alt(q_{2(n-k),\bi}) ,
\end{equation}
where the multilinear functionals $p_{2k,a}, q_{2k,a} \in \otimes^{2k} (\R^{16})^\ast$ are given by
\begin{align*}
p_{2k,a}(X_1, \dots, X_{2k}) & \coloneqq \bigl\langle P_{2k}(x_{I_1}, \dots, x_{I_{2k}}), a \bigr\rangle, \\
q_{2k,a}(X_1, \dots, X_{2k}) & \coloneqq \bigl\langle P_{2k}(y_{I_1}, \dots, y_{I_{2k}}), a \bigr\rangle,
\end{align*}
for $a = \bone,\bi$.
\end{Lem}

\begin{proof}
Referring back to \eqref{eq:clifford8model}, we note that for $u,v \in \bO$, the following holds:
\[
\begin{pmatrix}
0 & R_u \\
-R_{\bar{u}} & 0
\end{pmatrix}
\begin{pmatrix}
	0 & R_v \\
	-R_{\bar{v}} & 0
\end{pmatrix}
=
\begin{pmatrix}
-R_u R_{\bar{v}} & 0 \\
0 & -R_{\bar{u}}R_v
\end{pmatrix} .
\]
When applied to a positive spinor $w \in \bO^+$, only the top-left block is relevant, yielding
\begin{equation}
	\label{eq:octonion_product}
-R_u R_{\bar{v}}w = -(w \bar{v})u .
\end{equation}
Given the definitions of $s$ and $t$ from the beginning of the section, the multilinear functional
\[
T_{2n}(X_1, \dots, X_{2n}) \coloneqq \langle X_1 \cdots X_{2n} s, t \rangle,
\]
from Lemma~\ref{lem:clifford_prod}, is equal to
\begin{align*}
& \sum_{a \in \{\bone, \bi\}} \sum_{k=0}^n \sum_{I \in \Lambda(2n,2k)} \sgn(\sigma_I) \bigl\langle P_{2k}(x_{I_1}, \dots, x_{I_{2k}}), a \bigr\rangle \bigl\langle P_{2(n-k)}(y_{I_1^c}, \dots, y_{I_{2(n-k)}^c}), a \bigr\rangle .
\end{align*}
Since $n$ is even, the negative signs from \eqref{eq:octonion_product} cancel out in pairs. Applying Lemma~\ref{lem:tensor_alternation} to alternate the expression for $T_{2n}$ derived above, we obtain
\begin{align}
	\label{eq:alternation}
\Alt(T_{2n}) = \sum_{k=0}^n \Alt(p_{2k,\bone}) \wedge \Alt(q_{2(n-k),\bone}) + \Alt(p_{2k,\bi}) \wedge \Alt(q_{2(n-k),\bi}) .
\end{align}
If $X_1, \dots, X_{2n} \in \R^{16}$ are orthonormal, then $X_i X_j = -X_j X_i$ for different $i,j$, and therefore
\begin{align*}
T_{2n}(X_1, \dots, X_{2n}) & = \frac{1}{(2n)!} \sum_{\sigma \in S_{2n}} \sgn(\sigma) \langle X_{\sigma(1)} \cdots X_{\sigma(2n)} s, t \rangle \\
& = \Alt(T_{2n})(X_1, \dots, X_{2n}) .
\end{align*}
Using \eqref{eq:phi_2n} and \eqref{eq:alternation}, we obtain \eqref{eq:phi_identity} as stated, since
\[
\phi_{2n}(X_1 \wedges X_{2n}) = \langle X_1\cdots X_{2n} s, t\rangle = \Alt(T_{2n})(X_1, \dots, X_{2n}).
\]
Initially, this holds when evaluated on orthonormal systems. However, since both sides of \eqref{eq:phi_identity} are alternating, it extends naturally to all systems.
\end{proof}

The next step is to determine these alternations explicitly. We begin with $\phi_4$. By symmetry, it suffices to consider the tensors in the first factor of $\R^{16}$. By definition,
\begin{align*}
p_{0,a} & = \langle \bone, a\rangle , \\
p_{2,a}(x_1,x_2) & = \langle \bar{x}_2x_1, a \rangle , \\
p_{4,a}(x_1,x_2,x_3,x_4) & = \langle ((\bar{x}_4 x_3)\bar{x}_2)x_1, a \rangle .
\end{align*}
Clearly, $p_{0,\bi} = 0$ and for $p_{2,\bone}$ we have
\[
p_{2,\bone}(x_1,x_2) = \langle \bar{x}_2x_1, \bone \rangle = \langle x_1, x_2 \rangle = \langle \bone, \bar{x}_1x_2 \rangle = p_{2,\bone}(x_2,x_1) .
\]
Since $p_{2,\bone}$ is symmetric, its alternation vanishes. As in Section~\ref{subsec:classical_calibrations}, let $J$ be the complex structure on $\bO$, defined by $Jv \coloneqq v\bi$ with corresponding Kähler form $\omega \in \ext^2 \R^8$ given by
\[
\omega(v,w) \coloneqq \langle Jv, w\rangle = \langle v \bi, w \rangle = \langle \bi, \bar{v}w \rangle = p_{2,\bi}(w,v) .
\]
Since $p_{2,\bi}$ is already alternating, we have
\[
\Alt(p_{2,\bi}) = p_{2,\bi} = -\omega .
\]
Using \eqref{eq:algebra_basics} and \eqref{eq:cross}, we compute
\begin{align*}
p_{4,\bone}(x_1,x_2,x_3,x_4) & = \langle ((\bar{x}_4 x_3)\bar{x}_2)x_1 , \bone \rangle = \bigl\langle \overline{((\bar{x}_4 x_3)\bar{x}_2)x_1}, \bone \bigr\rangle = \langle \bar{x}_1(x_2(\bar{x}_3 x_4)), \bone \rangle .
\end{align*}
By Lemma~\ref{lem:cross_alternation}, the alternation $\Alt(p_{4,\bone})$ is the Cayley calibration $\phi_{\mathrm{Cay}} \in \ext^4\R^8$. The holomorphic volume form is given by $\Omega \coloneqq Z_{1234} \in \ext^4\R^8$ in terms of complex coordinates as in \eqref{eq:complex_coordinates}. For $V_8 \coloneqq \R^8 \times \{0\}$ this yields the forms $\Omega$, $\omega$ and $\phi_{\mathrm{Cay}}$ related to the complex structure $J$. The corresponding forms on $V_8' \coloneqq \{0\} \times \R^8$ are denoted by $\Omega'$, $\omega'$ and $\phi_{\mathrm{Cay}}'$ with complex structure $J'$. Using \eqref{eq:phi_identity} and \eqref{eq:cayley_identity}, we obtain
\begin{align}
	\nonumber
	\phi_4 & = \phi_{\mathrm{Cay}} + \phi_{\mathrm{Cay}}' + \omega \wedge \omega' \\
	\nonumber
	& = \Re \Omega - \tfrac{1}{2}\omega^2 + \Re \Omega' - \tfrac{1}{2}\omega'^2 + \omega \wedge \omega' \\
	\label{eq:phi4}
	& = \Re \Omega + \Re \Omega' - \tfrac{1}{2}(\omega - \omega')^2.
\end{align}

Now, we proceed with $\phi_8$. We already derived $p_{0,\bi} = q_{0,\bi} = 0$ and $\Alt(p_{2,\bone}) = \Alt(q_{2,\bone}) = 0$. Applying \eqref{eq:phi_identity} as before, we obtain
\begin{align*}
\phi_8 & = \sum_{k=0}^4 \Alt(p_{2k,\bone}) \wedge \Alt(q_{2(n-k),\bone}) + \Alt(p_{2k,\bi}) \wedge \Alt(q_{2(n-k),\bi}) \\
 & = \Alt(p_{8,\bone}) + \Alt(q_{8,\bone}) + \Alt(p_{4,\bone}) \wedge \Alt(q_{4,\bone}) + \Alt(p_{4,\bi}) \wedge \Alt(q_{4,\bi}) \\
 & \phantom{=}\ + \Alt(p_{6,\bi}) \wedge \Alt(q_{2,\bi}) + \Alt(p_{2,\bi}) \wedge \Alt(q_{6,\bi}) .
\end{align*}
We have already determined that $\Alt(p_{2,\bi}) = -\omega$ and $\Alt(p_{4,\bone}) = \phi_{\mathrm{Cay}}$. The remaining forms are identified in the following lemma.

\begin{Lem}
$\Alt(p_{8,\bone}) = \vol = \frac{1}{4!}\omega^4$, $\Alt(p_{6,\bi}) = \frac{1}{3!}\omega^3$ and $\Alt(p_{4,\bi}) = \Im\Omega$.
\end{Lem}

\begin{proof}
First, observe that each $P_{2k}$ and $p_{2k,a}$ is already alternating when evaluated on orthogonal systems, due to \eqref{eq:oct_permute}. We have
\begin{align*}
p_{4,\bi}(x_1,x_2,x_3,x_4) & = \langle ((\bar{x}_4 x_3)\bar{x}_2)x_1 , \bi \rangle = \bigl\langle \overline{((\bar{x}_4 x_3)\bar{x}_2)x_1}, \bar{\bi} \bigr\rangle \\
& = -\langle \bar{x}_1(x_2(\bar{x}_3 x_4)), \bi \rangle .
\end{align*}
Thus $\Alt(p_{4,\bi}) = \Im\Omega$ as a consequence of Lemma~\ref{lem:cross_i}, and
$\Alt(p_{8,\bone}) = \vol$ follows from the fact that $P_8(e_1, \dots, e_8) = \bone$, as established in Lemma~\ref{lem:octonion_product}. The identity $\Alt(p_{6,\bi}) = \frac{1}{3!}\omega^3$ follows by duality. For $I \in \Lambda(8,6)$, the above evaluation of $P_8$ yields
\[
\sgn(\sigma_I) \bone = P_8(e_{I_1}, \dots, e_{I_6},e_{I^c_1},e_{I^c_2}) = (P_6(e_{I_1}, \dots, e_{I_6})\bar{e}_{I^c_2})e_{I^c_1}.
\]
Using \eqref{eq:algebra_basics}, we find
\begin{equation*}
\sgn(\sigma_I) p_{2,\bi}(e_{I_2^c},e_{I_1^c}) = \sgn(\sigma_I) \bigl\langle \bar{e}_{I_1^c} e_{I_2^c}, \bi \bigr\rangle = p_{6,\bi}(e_{I_1}, \dots, e_{I_6}) .
\end{equation*}
Since $p_{2,\bi} = -\omega$, it follows from Lemma~\ref{lem:standard_calibrations} that $\Alt(p_{6,\bi}) = \ast \omega = \tfrac{1}{3!}\omega^3$.
\end{proof}

By collecting all the terms, we can now compute $\phi_8$ as follows:
\begin{align}
	\nonumber
\phi_8 & = \tfrac{1}{4!}\omega^4 + \tfrac{1}{4!}\omega'^4 - \tfrac{1}{3!} \omega \wedge \omega'^3 - \tfrac{1}{3!} \omega^3 \wedge \omega' \\
	\nonumber
& \phantom{=}\ + \Im\Omega \wedge \Im\Omega' + \phi_{\mathrm{Cay}} \wedge \phi_{\mathrm{Cay}}' \\
	\nonumber
& = \tfrac{1}{4!}(\omega^4 - 4\, \omega \wedge \omega'^3 + 6\, \omega^2 \wedge \omega'^2 - 4\, \omega^3 \wedge \omega' + \omega'^4) \\
	\nonumber
& \phantom{=}\ + \Im \Omega \wedge \Im\Omega' + \Re \Omega  \wedge \Re \Omega'  \\
	\nonumber
& \phantom{=}\ - \Re\Omega \wedge \tfrac{1}{2}\omega'^2 - \tfrac{1}{2}\omega^2 \wedge \Re\Omega' \\
\label{eq:phi8}
& = \tfrac{1}{4!}(\omega-\omega')^4 + \Im\Omega \wedge \Im\Omega' + \Re\Omega \wedge \Re\Omega' \\
\nonumber
& \phantom{=}\ - \Re\Omega \wedge \tfrac{1}{2}\omega'^2 - \tfrac{1}{2}\omega^2 \wedge \Re\Omega' .
\end{align}
Gathering all this information, we arrive at a proof of  Theorem~\ref{thm:spinorproduct}.

\begin{proof}[Proof of Theorem~\ref{thm:spinorproduct}]
The element $\phi \in \ext^\ast \R^{16}$ decomposes into $\psi \coloneqq 16^2\, s\scirc s$ and $\eta \coloneqq 16^2\, s'\scirc s$. Since $p_{0,\bi} = q_{0,\bi} = 0$, we obtain $\eta_0 = \eta_{16} = 0$. The following properties of $\phi$ are taken from Theorem~\ref{thm:spinorprod_dh}:
\begin{itemize}
	\item $\psi_k$ and $\eta_k$ are calibrations for all $k=0, \dots, 16$.
	\item $\phi_k = \psi_k = \eta_k = 0$ in case $k$ is odd.
	\item $\phi_k = \eta_k$ in case $k \equiv 2 \text{ mod } 4$.
	\item $\phi_0 = \psi_0 = 1$, $\phi_{16} = \psi_{16} = \lambda$, the volume form on $\R^{16}$.
	\item $\phi_{12} = \ast \phi_4$.
\end{itemize}
Thus, in order to show that $\phi_k$ is a calibration for every $k$, it remains to check it for $k=4,8$.

Reversing the complex structure on $V_8'$ from $J'$ to $-J'$ leaves $\Re\Omega'$ unchanged, while $\omega'$ and $\Im\Omega'$ change sign. The form $\phi_4$, as given in \eqref{eq:phi4}, is thus
\[
\phi_4 = \Re\Omega + \Re\Omega' - \tfrac{1}{2}(\omega + \omega')^2
\]
with respect to the complex structure $\bar{J} \coloneqq (J,-J')$ on $\R^{16}$. It is straightforward to verify that $|\phi_4(\xi)| \leq 1$ for $\xi \in \bG(\R^{16},4)$ by expressing $\xi$ in normal form with respect to the complex structure $\bar{J}$. Moreover, $\xi$ is calibrated by $\phi_4$ if and only if $\xi$ is a Cayley $4$-plane in $V_8$ or $V_8'$; or a negatively oriented complex $2$-plane with respect to $\bar{J}$. The argument follows similarly to that of \eqref{eq:extreme_case_estimate} and is left to the reader. 

The form $\phi_8$ is computed in \eqref{eq:phi8}. With respect to the modified complex structure $\bar{J}$ on $\R^{16}$, we obtain
\begin{align*}
\phi_8 & = \tfrac{1}{4!}(\omega + \omega')^4 + \Re\Omega \wedge \Re\Omega' - \Im\Omega \wedge \Im\Omega' \\
& \phantom{=}\ - \Re\Omega \wedge \tfrac{1}{2}\omega'^2 - \tfrac{1}{2}\omega^2 \wedge \Re\Omega' \\
& = \Re(\Omega \wedge \Omega') + \tfrac{1}{4!}(\omega + \omega')^4 - \Re\Omega \wedge \tfrac{1}{2}\omega'^2 - \tfrac{1}{2}\omega^2 \wedge \Re\Omega' .
\end{align*}
If, in the identification of $(\R^{16},\bar{J})$ with $\C^8$, the coordinates $z_1$ and $z_5$ of $\C^8$ are multiplied by $-1$, then $\Omega$ and $\Omega'$ change signs, while all other forms remain unchanged. As a result,
\[
\phi_8 = \Re(\Omega \wedge \Omega') + \tfrac{1}{4!}(\omega + \omega')^4 + \Re\Omega \wedge \tfrac{1}{2}\omega'^2 + \tfrac{1}{2}\omega^2 \wedge \Re\Omega'
\]
coincides with $\Phi$ in Theorem~\ref{thm:main}, confirming that $\phi_8$ is a calibration.
\end{proof}

\begin{Rem}
It may be insightful to compare the contributions to $|\phi_k|^2$ from the two components $\psi_k$ and $\eta_k$. According to \cite[Proposition~2.1(7)]{HL}, we have $|\psi|^2 = |\eta|^2 = 16^2$. Since elements of different degree in $\Cl(16)$ are orthogonal, it follows that
\[
\sum_k |\psi_k|^2 = \sum_k |\eta_k|^2 = 256 .
\]
Expanding $\psi$ and $\eta$ in the standard basis of $\ext^\ast\R^{16}$, as in the proof of Theorem~\ref{thm:main}, yields the following squared norms in even degrees $k=0,2, \dots, 16$:

\begin{table}[!ht]
\centering
\begin{tabular}{c|ccccccccc|c}
\hline
~ & $0$ & $2$ & $4$ & $6$ & $8$ & $10$ & $12$ & $14$ & $16$ & $\sum$ \\ \hline
$\psi$ & 1 & 0 & 28 & 0 & 198 & 0 & 28 & 0 & 1 & 256 \\ 
$\eta$ & 0 & 0 & 16 & 64 & 96 & 64 & 16 & 0 & 0 & 256 \\ 
$\phi$ & 1 & 0 & 44 & 64 & 294 & 64 & 44 & 0 & 1 & 512 \\ \hline
\end{tabular}
\end{table}

The squared norms $|\psi_4|^2 = 28$ and $|\eta_4|^2 = 16$ follow from \eqref{eq:phi4}, noting that
\[
\psi_4 = \phi_{\mathrm{Cay}} + \phi_{\mathrm{Cay}}' \quad \text{and} \quad \eta_4 = \omega \wedge \omega'.
\]
The values $|\psi_8|^2 = 198$ and $|\eta_8|^2 = 96$ are obtained by identifying
\[
\psi_8 = \tfrac{1}{4!}\omega^4 + \tfrac{1}{4!}\omega'^4 + \phi_{\mathrm{Cay}} \wedge \phi_{\mathrm{Cay}}',
\]
\[
\eta_8 = \Im\Omega \wedge \Im\Omega' - \omega \wedge \tfrac{1}{3!} \omega'^3 - \tfrac{1}{3!} \omega^3 \wedge \omega',
\]
in \eqref{eq:phi8}. Furthermore, since $p_{0,\bi} = q_{0,\bi} = 0$, it follows that $\eta_2 = 0$. Given the duality relation $\ast \psi_{k} = -\psi_{16-k}$ for $k \equiv 2 \text{ mod } 4$ and the fact that $|\psi|^2 = |\eta|^2 = 256$, the remaining norms must satisfy $|\eta_6|^2 = |\eta_{10}|^2 = 64$.
\end{Rem}

\vspace{0.0cm}

\subsection*{Acknowledgment}
I thank the referees for their careful reading and helpful comments. Special thanks go to the referee who provided an independent verification of Theorem 1.1.

%%%%%%%%%%%%%%%%%%%%%%%%%%%%%%%%%%%%%%%%%%%%%%%%%%%%%%%%%%%%%%%%%%%%%%%%%%%%%%%

\end{document}